\documentclass[11pt]{amsart}
\usepackage{amscd,amssymb,amsxtra}
\usepackage[mathscr]{eucal}
\usepackage{comment}
\usepackage{booktabs}
\makeatletter
\def\section{\@startsection{section}{1}%
  \z@{4\linespacing\@plus\linespacing}{\linespacing}%
  {\normalfont\scshape\centering}}
 \def\subsection{\@startsection{subsection}{2}%
   \z@{1.25\linespacing\@plus.7\linespacing}{.5\linespacing}%
   {\normalfont\bfseries}}
\makeatother

\theoremstyle{definition}
\newtheorem{example}[equation]{Example}

\newtheorem*{example*}{Example}
\newtheorem*{problem*}{Problem}
\newtheorem*{exercise*}{Exercise}
\newtheorem*{question*}{Question}

\theoremstyle{remark}

\newtheorem{remark}[equation]{Remark}

\newtheorem*{note*}{Note}
\newtheorem*{notation*}{Notation}
\newtheorem*{claim*}{Claim}
\newtheorem*{remark*}{Remark}

\theoremstyle{plain}
\newtheorem{definition}[equation]{Definition}
\newtheorem{theorem}[equation]{Theorem}

\newtheorem{proposition}[equation]{Proposition}
\newtheorem{conjecture}[equation]{Conjecture}
 
\newtheorem*{definition*}{Definition}
\newtheorem*{theorem*}{Theorem}
\newtheorem*{corollary*}{Corollary}
\newtheorem*{lemma*}{Lemma}
\newtheorem*{proposition*}{Proposition}
\newtheorem*{conjecture*}{Conjecture}

\numberwithin{equation}{section}

\renewcommand{\:}{\colon}

\newcommand{\CC}{{\mathbb C}}
\newcommand{\CP}{{\mathbb C\mathbb P}}

\DeclareMathOperator{\Diff}{Diff}

\newcommand{\FF}{\mathbb F}

\DeclareMathOperator{\Hom}{Hom}
\DeclareMathOperator{\id}{id}

\DeclareMathOperator{\pt}{pt}

\newcommand{\RR}{{\mathbb R}}

\newcommand{\ZZ}{{\mathbb Z}}

\newcommand{\chiup}{\raise.5ex\hbox{$\chi$}}
\newcommand{\cir}{S^1}

\newcommand{\inv}{^{-1}}
\newcommand{\mstrut}{^{\vphantom{1*\prime y}}}

\DeclareMathOperator{\rank}{rank}
\newcommand{\res}[1]{\negmedspace\bigm|\mstrut_{#1}}
\newcommand{\temsquare}{\raise3.5pt\hbox{\boxed{ }}}

\newcommand{\zmod}[1]{\ZZ/#1\ZZ}
\newcommand{\zn}{\zmod{n}}
\newcommand{\zt}{\zmod2}

\usepackage[all,2cell,dvips]{xy}\renewcommand{\cir}{\ensuremath{S^1}}
\UseAllTwocells
\usepackage{graphicx}
\usepackage{epsf}
\usepackage{pstricks}
\usepackage{pst-grad}

\DeclareMathOperator{\Ad}{Ad}
\DeclareMathOperator{\Lie}{Lie}
\DeclareMathOperator{\module}{-mod}
\DeclareMathOperator{\trace}{trace}
\DeclareMathOperator{\vol}{vol}
\newcommand{\AP}{\mathcal{A}\mstrut_P}
\newcommand{\ApP}{\mathcal{A}\mstrut_{\pi ^*P}}
\newcommand{\C}{\mathscr{C}}
\newcommand{\cat}{\mathcal{C}}
\newcommand{\FX}{\mathcal{F}\mstrut _X}
\newcommand{\FY}{\mathcal{F}\mstrut _Y}
\newcommand{\FbX}{\mathcal{F}\mstrut _{\bX}}
\newcommand{\Ftil}{\tilde{F}}
\newcommand{\Hilb}[1]{\mathcal{H}_{#1}}
\newcommand{\Ouniv}{\Omega \mstrut _P}
\newcommand{\R}{\mathscr{R}}
\newcommand{\TD}{\Theta \mstrut _\Delta }
\newcommand{\Tuniv}{\Theta \mstrut _P}
\newcommand{\Tun}{\Theta ^{\textnormal{univ}}}
\newcommand{\X}{\mathscr{X}}
\newcommand{\Y}{\mathscr{Y}}
\newcommand{\bM}{\partial \Sigma }
\newcommand{\bX}{\partial X}
\newcommand{\bg}{\bar{\gamma }}
\newcommand{\bordO}{\textnormal{Bord}^{O}} 
\newcommand{\bordSO}{\textnormal{Bord}^{SO}} 
\newcommand{\bordSpin}{\textnormal{Bord}^{Spin}} 
\newcommand{\bordfr}{\textnormal{Bord}^{\textnormal{fr}}} 
\newcommand{\fconn}[1]{\mathcal{M}_{#1}}
\newcommand{\field}[1]{\mathcal{F}_{#1}}
\newcommand{\gl}{(G,\lambda )}
\newcommand{\gpd}{/\!/}
\newcommand{\hY}{\hat{Y}}

\newcommand{\multc}{
      \pscustom[fillstyle=gradient, 
    gradbegin=white, gradend=gray,gradmidpoint=0,gradangle=70]{
        \psbezier(1.5,2.5)(1.5,1.1)(.4,1.6)(.5,0)
        \psbezier(.5,0)(.4,-.25)(-.4,-.25)(-.5,0)
        \psbezier(-0.5,0)(-.4,1.6)(-1.5,1.1)(-1.5,2.5)
        \psline(-.5,2.5)
        \psbezier(-.5,2.5)(-.6,1.5)(0.6,1.5)(.5,2.5)
        \psline(1.5,2.5)
    }
    \psellipse[fillcolor=lightgray,fillstyle=gradient,
        gradbegin=lightgray, gradend=gray,gradmidpoint=1,gradangle=110](-1,2.5)(.5,.2)
    \psellipse[fillcolor=lightgray,fillstyle=gradient,
        gradbegin=lightgray, gradend=gray,gradmidpoint=1,gradangle=110](1,2.5)(.5,.2)
     \begin{psclip}{
 \pspolygon[linestyle=none](.5,0)(.5,.3)(-.5,.3)(-.5,0)(.5,0)
 }
 \psellipse[linestyle=dotted](0,0)(.5,0.2)
 \end{psclip}
 }
 
\newcommand{\birthc}{
 \pscustom[fillstyle=gradient,
    gradbegin=white, gradend=gray,gradmidpoint=0,gradangle=110]{
        \psbezier(-.5,0)(-.5,.9)(0.5,.9)(.5,0)
        \psbezier(.5,0)(.4,-.25)(-.4,-.25)(-.5,0)
    }
 \begin{psclip}{
 \pspolygon[linestyle=none](.5,0)(.5,.3)(-.5,.3)(-.5,0)(.5,0)
 }
 \psellipse[linestyle=dotted](0,0)(.5,0.2)
 \end{psclip}
 }

\theoremstyle{plain}
\newtheorem*{observation*}{Observation}

\begin{document}

\abovedisplayskip18pt plus4.5pt minus9pt
\belowdisplayskip \abovedisplayskip
\abovedisplayshortskip0pt plus4.5pt
\belowdisplayshortskip10.5pt plus4.5pt minus6pt
\baselineskip=15 truept
\marginparwidth=55pt




 \title[Chern-Simons Theory]{Remarks on Chern-Simons Theory}
 \author[D. S. Freed]{Daniel S.~Freed}
 \thanks{The author is supported by NSF grant DMS-0603964}
 \address{Department of Mathematics \\ University of Texas \\ 1 University
Station C1200\\ Austin, TX 78712-0257}
 \email{dafr@math.utexas.edu}
 \dedicatory{Dedicated to MSRI on its $25^{\text{th}}$ anniversary}
 \date{August 8, 2008; revised October 26, 2008}
 \thanks{Based on a talk given in the Simons Auditorium in Chern Hall at the
Mathematical Sciences Research Institute on the occasion of its
$25^{\text{th}}$~anniversary.}
\maketitle

 
A beautiful line of development in Riemannian geometry is the relationship
between curvature and topology.  In one of his first major works, written
in~1946, Chern~\cite{Ch} proves a generalized Gauss-Bonnet theorem by
producing what we now call a transgressing form on the unit sphere bundle of
the manifold.  Twenty-five years later, together with Simons, he took up
transgression~\cite{CSi} in the context of the theory of connections on
arbitrary principal bundles.  The Chern-Simons invariants of a connection are
secondary \emph{geometric} invariants; in between these works Chern and Weil
developed the theory of primary \emph{topological} invariants of connections.
Both invariants are local in the sense that they are computed by integrals of
differential forms.  The relationship between them is that the differential
of the Chern-Simons form is the Chern-Weil form.  The integral of the
Chern-Weil form over a closed manifold is independent of the connection, so
is a topological invariant.  In the late~1980s Edward Witten~\cite{W1}
proposed a new topological invariant of 3-manifolds from these same
ingredients.  He achieves topological invariance via a technique unavailable
to geometers: Witten integrates the exponentiated Chern-Simons invariant over
the \emph{infinite-dimensional} space of all connections.  Because the
connection is integrated out, the result depends only on the underlying
manifold.  This, then, is the quantum Chern-Simons invariant.
 
There are rich stories to tell about both the classical and quantum
Chern-Simons invariants in geometry, topology, and physics.  The classical
Chern-Simons invariant is an obstruction to conformal immersions of
3-manifolds into Euclidean space, is closely related to the
Atiyah-Patodi-Singer invariant, and was refined in the Cheeger-Simons theory
of differential characters.  It appeared in physics before Witten's work, for
example in the theory of anomalies.  The quantum Chern-Simons invariant is
closely related to the Jones invariants~\cite{Jo1} of links, which have had
many applications in knot theory.  But we do not attempt here a review of all
work on Chern-Simons invariants.  Our interest is the \emph{structure} behind
the quantum invariants as they relate to structure in physics.  Indeed,
Witten's achievement was to fit the Jones invariants into a larger structure,
that of a 3-dimensional \emph{quantum field theory}.
 
The beginning of the $20^{\text{th}}$~century saw two revolutionary
developments in physics: relativity and quantum mechanics.  The mathematical
basis for general relativity---differential geometry---was already
well-developed when Einstein came along, and in turn his gravitational
equations, which focus attention on Ricci curvature, spurred many
developments in Riemannian geometry.  Quantum mechanics, on the other hand,
quickly inspired the development of operator theory and parts of
representation theory.  These foundations for quantum mechanics later
influenced many other parts of mathematics, including diverse areas such as
partial differential equations and number theory~\cite{Ma}.  Quantum field
theory in its first incarnation, quantum electrodynamics, combines quantum
mechanics and Maxwell's classical theory of the electromagnetic field.  It
has provided some of the most precise computational agreement between theory
and experiment, and quantum field theory is the setting for the standard
model of particle physics.  Although there has been much mathematical work on
quantum field theory, its foundations are not at all settled.  The
interaction with mathematics has greatly broadened over the past 25~years:
now quantum field theory---and also string theory---enjoys a deep interaction
with many branches of mathematics, suggesting novel results, surprising
connections, and new lines of research.  One important impetus for this
development was the advent of new examples of quantum field theories closely
connected with geometry and topology.  Quantum Chern-Simons theory is one of
the first examples, and it is a purely topological theory at both the
classical and quantum levels.  We use it here as a focal point to discuss
topological quantum field theories (TQFT) in general.  From it we have gained
insight into the formal structure of all quantum field theories, not just
topological ones, and it is that general structure which we accentuate.
Also, Chern-Simons theory provides a window into the entire interchange
between mathematics and these parts of physics.  We return to this line of
thought in the second half of~\S\ref{sec:7}.
 
Witten's main tool is the \emph{path integral}, an integral on a function
space.  (In Feynman's approach to quantum mechanics~\cite{Fe} it is truly an
integral over a space of paths; in general field theories it is still usually
referred to as the `path integral', though `functional integral' is better.)
Regrettably, the measures needed for such integrals have only been
constructed in special theories, and these do not cover the example of
Chern-Simons theory.  In~\S\ref{sec:5} we extract some important formal
properties of path integrals in general.  We focus on \emph{locality}, which
is manifested in gluing laws.  This formal structure is nicely expressed in
familiar mathematical terms, as a \emph{linearization of correspondence
diagrams}.  Here the correspondence diagrams are built from the semiclassical
fields; the linearization is the quantum theory.  Following Segal~\cite{Se2}
and Atiyah~\cite{A2} we abstract axioms for a TQFT, which we give
in~\S\ref{sec:1}.  The axioms do not encode the passage from fields to the
quantum theory, but only the properties of the quantum theory itself.  In the
end these axioms are quite analogous to those of homology theory, but with
important differences: (i)\ a TQFT is defined on manifolds of a particular
dimension, whereas homology theory is defined on arbitrary topological
spaces; and (ii)\ a TQFT is multiplicative (say on disjoint unions) whereas
homology theory is additive.  The axioms are most neatly stated in
categorical language, which, given the analogy with homology theory, is not a
surprise.
 
While the definition of a TQFT is simple, examples are not so easy to come
by.  We distinguish between \emph{generators-and-relations} constructions and
\emph{a priori} constructions.  The former are most prevalent in the
mathematics literature, and we present some theorems which classify TQFTs in
these terms.  (See Theorem~\ref{thm:3} and Theorem~\ref{thm:5}, for example.)
Most rigorous constructions of the 3-dimensional Chern-Simons theory are of
this generators-and-relations type.  The path integral is an \emph{a priori}
construction---the axioms of a TQFT follow directly and geometrically--- but
in most cases it is not rigorous; a notable exception is gauge theory with
\emph{finite} gauge group where the path integral reduces to a finite sum.
The 2-dimensional reduction of Chern-Simons theory for \emph{any} gauge group
has an \emph{a priori} construction.  In fact, that reduced theory exists
over the integers, whereas the 3-dimensional Chern-Simons theory is defined
over the complex numbers.  It may be that TQFTs which are dimensional
reductions always exhibit integrality---and even a connection to
$K$-theory---as we speculate in~\S\ref{sec:6}.
 
Another theme here is the extension of a TQFT to higher codimension.  Usually
in a quantum field theory there are two dimensions in play: the dimension~$n$
of spacetime and the dimension~$n-1$ of space.  For sure locality in space
appears in the physics literature---think of the Hilbert space of a lattice
model constructed as a tensor product of local Hilbert spaces, or of local
algebras of observables in some axiomatic treatments of quantum field
theory---though that aspect is not always emphasized.  One of the very
fruitful idea in \emph{topological} theories is to consider gluing laws along
corners.  In the best case one goes all the way down to points.  So, for
example, the usual 3-dimensional Chern-Simons theory may be called a
`2-3~theory' as it involves 2-~and 3-dimensional manifolds.  The extension
down to points is called a `0-1-2-3~theory'.  Theories which extend down to
points truly earn the adjective `local'; other theories are only partially
local.  In this connection Segal~\cite{Se2} has observed that the existence
of handlebody decompositions for manifolds means that in a theory with at
least three ``tiers'' all of the data is determined in principle from balls
and products of spheres.  In another direction one can often build into a
topological quantum field theory invariants for families of manifolds as well
as invariants of a single manifold.  String topology as well as the theory
which encodes Gromov-Witten invariants both include invariants of families.
One new development in this area are generators-and-relations structure
theorems for TQFTs which are fully local---go down to points---and include
families: see Theorem~\ref{thm:29}, Theorem~\ref{thm:31}, and
Conjecture~\ref{thm:30}.  In these statements the values of a TQFT are not
necessarily Hilbert spaces and complex numbers.  For these are really
generators-and-relations results about bordism categories of manifolds, so
the codomain can be very general.  This flexibility in the codomain is
necessary if one is to avoid semisimplicity when considering theories over
rings which are not fields (see Remark~\ref{thm:2}).
 
In an $n$-dimensional theory the invariant of an $n$-manifold is a
\emph{number} and the invariant of an $(n-1)$-manifold is a \emph{set} (the
quantum Hilbert space).  It is natural, then, that the invariant of an
$(n-2)$-manifold be a \emph{category}; the linear structure of quantum theory
demands that it be a \emph{linear} category.  Further descent in dimension,
obligatory for full locality, requires a concomitant ascent in category
number.  This is a familiar ladder one climbs---or descends, depending on
your orientation---in topology.  The simplest invariant of a space~$S$ is the
\emph{set}~$\pi _0S$ of path components.  The next invariant is a
\emph{category}: the fundamental groupoid~$\pi _{\le 1}S$.  At the next level
we encounter a \emph{2-category}~$\pi _{\le2}S$ and so on.  To incorporate
families of manifolds one enriches to topological categories or some
equivalent.

Chern-Simons theory as initially conceived from the path integral is a
2-3~theory.  Reshetikhin and Turaev~\cite{RT}, \cite{Tu} give a rigorous
mathematical construction when the gauge group is simple and simply
connected; it utilizes a suitable category of representations of a related
quantum group.  Hence they construct a 1-2-3~theory; this category is
attached to~$\cir$.  It has been a longstanding question to extend to a
0-1-2-3~theory: What does Chern-Simons attach to a point?  In~\S\ref{sec:2}
we give the answer for a finite gauge group (it was already contained
in~\cite{F2}).  There is an extension of these ideas which works when the
gauge group is a torus group, and hopefully in more general cases as well.
 
In Chern-Simons theory the vector space attached to the torus~$\cir\times
\cir$ refines to a free abelian group which carries a ring structure, the
so-called \emph{Verlinde ring}.  In a series of papers~\cite{FHT1} this ring
was identified with a certain $K$-theory ring built directly from the gauge
group; see Theorem~\ref{thm:18}.  This opens the possibility of an \emph{a
priori} construction of quantum Chern-Simons theory using pure topology,
specifically $K$-theory.  There is an \emph{a priori} construction of the
2-dimensional reduction, as recounted at the end of~\S\ref{sec:6}, and it
\emph{is} pure $K$-theory.  It remains to be seen if this can be extended to
an \emph{a priori} 3-dimensional construction, or at least if the input to
the generators-and-relations Conjecture~\ref{thm:30} can be given in terms of
$K$-theory (as it is for finite groups and tori).
 
We begin in~\S\ref{sec:3} with a brief account of the differential geometry
into which the classical Chern-Simons invariant fits; the actual invariant is
constructed in the appendix.  The mathematics which goes into the rigorous
constructions of the quantum theory---quantum groups and the like---are quite
different: only the path integral approach brings in the differential
geometry directly.  We emphasize this in~\S\ref{sec:7} where we observe that
an elementary consequence of the path integral~\eqref{eq:58}, which involves
classical invariants of flat connections, is not captured by the axiomatics
and in fact remains unproven.  This state of affairs returns us to a more
general discussion about the interaction of mathematics with quantum field
theory and string theory.  The casual reader may wish to rejoin us at that
point.  In between we extract the formal structure of path
integrals~(\S\ref{sec:5}), present axiomatics for TQFT together with
generators-and-relations theorems~(\S\ref{sec:1}), discuss constructions of
Chern-Simons as a 1-2-3 and 0-1-2-3~theory~(\S\ref{sec:2}), and present the
relationship to $K$-theory~(\S\ref{sec:6}).
 
I had the good fortune to be a graduate student at Berkeley when MSRI opened
its doors in~1982.  From the beginning Chern nurtured an open and welcoming
atmosphere.  I quickly discovered that whereas faculty members at Evans
Hall---as at any mathematics department---were occupied with teaching and
administration as well as research, the visitors to MSRI were engaging with
mathematics and (even young) mathematicians all day long.  I learned early
the world of difference between a mathematics department and a Mathematical
Sciences Research Institute.  No wonder that so many institutes around the
world emulate MSRI!  Over its first 25~years MSRI has grown tremendously in
scope and outreach all the while cultivating a creative environment for
mathematics and mathematicians to flourish.  I am very grateful for the time
I have been able to spend at MSRI, and it is an honor to dedicate this paper
to the continued good health of MSRI.
 
The circle of ideas surrounding Chern-Simons theory has been the topic of
discussions with my collaborators Michael Hopkins and Constantin Teleman for
at least eight years, and their influence can be felt throughout as can the
influence of earlier collaborators and my students.  For the more recent
topics Jacob Lurie and David Ben-Zvi have joined the conversation.  It is a
pleasure to thank them all.

  \section{The classical chern-simons invariant}\label{sec:3}

We begin with the $19^{\textnormal{th}}$~century progenitors of the
Chern-Weil and Chern-Simons work.  Let $\Sigma $~be a surface embedded in
3-dimensional Euclidean space.  The curvature of~$\Sigma $ is measured by a
single function $K\:\Sigma \to\RR$.  For example, the curvature of a sphere
of radius~$R$ is the constant function with value~$1/R^2$.  The curvature was
first investigated by Gauss in~1825.  His famous \emph{theorema
egregium}~\cite[p.~105]{G} proves that this \emph{Gauss curvature} is
intrinsic: it only depends on the induced metric on the surface, not on the
embedding of the surface into space.  The Riemannian metric determines a
measure~$d\mu\mstrut_\Sigma $ on~$\Sigma $.  The Gauss-Bonnet theorem,
apparently first given its global formulation by Walter van
Dyck~\cite[p.~141]{Hi}, states that if $\Sigma $~is \emph{closed}---compact
with no boundary---then
  \begin{equation}\label{eq:20}
     \int_{\Sigma }K\,d\mu\mstrut_\Sigma  =2\pi \chi (\Sigma ), 
  \end{equation}
where $\chi (\Sigma )$~is the Euler characteristic of~$\Sigma $.  This
formula is the apogee of a first course in differential geometry.  It and its
many generalizations link local geometry and global topology.

Now suppose $\Sigma $~has a boundary~$\partial \Sigma $ and is compact.  Then
there is an extra term which appears in~\eqref{eq:20}, the \emph{total
geodesic curvature} of the boundary.  Suppose $C\subset \Sigma $~is a closed
curve equipped with an orientation of its normal bundle.  The geodesic
curvature $\kappa \mstrut_C\:C\to\RR$ is a generalization of the curvature of
a plane curve; it vanishes if $C$~is a geodesic.  The metric induces a
measure~$d\mu \mstrut_{C}$ on~$C$ and
  \begin{equation}\label{eq:21}
     \int_{C}\kappa \mstrut_C\,d\mu \mstrut_C 
  \end{equation}
is the total geodesic curvature of~$C$.  The generalized Gauss-Bonnet formula
is  
  \begin{equation*}
     \int_{\Sigma }K\,d\mu \mstrut_\Sigma  + \int\mstrut_{\bM}\kappa \mstrut_{\bM}\,d\mu
     \mstrut_{\bM} = 2\pi \chi 
     (\Sigma ).  
  \end{equation*}

The classical Chern-Simons invariant~\cite{CSi} is a generalization of the
total geodesic curvature~\eqref{eq:21}.  The version we need is defined on a
compact oriented 3-manifold~$X$.  Like the total geodesic curvature it is an
\emph{extrinsic} invariant, but now the extrinsic geometry is defined by a
principal bundle with connection on~$X$, not by an embedding into a
Riemannian manifold.  The differential geometry is quite pretty, and for the
convenience of the reader we give a lightening review of connections and the
general Chern-Simons construction in the appendix.  For the exposition in the
next section we use a simplified version of this \emph{classical} invariant
as our emphasis is on its \emph{quantization} and ultimately the very
different mathematics used to compute in the quantum theory.

  \section{path integrals}\label{sec:5}
 Let $G=SU(n)$~be the Lie group of unitary $n\times n$~matrices of
determinant one for some~$n\ge2$.  Its Lie algebra~$\mathfrak{g}$ consists of
$n\times n$~skew-Hermitian matrices of trace zero.  Fix a closed oriented
3-manifold~$X$.  A connection~$A$ on the trivial $G$-bundle over~$X$ is a
skew-Hermitian matrix of 1-forms with trace zero, i.e., $A\in \Omega
^1(X;\mathfrak{g})$.  In this case the \emph{classical Chern-Simons
invariant} of~$A$ is given by the explicit formula
  \begin{equation}\label{eq:37}
     S(A) = \frac{1}{8\pi ^2}\int_{X}\trace\bigl( A\wedge dA \;+\; \frac 23
     \,A\wedge A\wedge A\bigr)
  \end{equation}
where the wedge products are combined with matrix multiplication.  The
integrand in~\eqref{eq:37} is a 3-form, and the integral depends on the
orientation of~$X$.  As mentioned earlier this is a 3-dimensional analog of
the 1-dimensional total geodesic curvature~\eqref{eq:21}; there are
generalizations to higher dimensions and all Lie groups, as explained in the
appendix.  Chern and Simons were particularly interested in the Levi-Civita
connection on the tangent bundle of a \emph{Riemannian} manifold~$X$, in
which case the Chern-Simons invariant is an obstruction to certain conformal
immersions.  They embarked on that study to derive combinatorial formulas for
the first Pontrjagin number of a compact oriented 4-manifold, or perhaps with
an eye toward the Poincar\'e conjecture.  Regardless, the classical
Chern-Simons invariant has found numerous applications in differential
geometry, global analysis, topology, and theoretical physics.
 
Edward Witten~\cite{W1} used the classical Chern-Simons invariant to derive
not an invariant of Riemannian 3-manifolds, but rather a topological
invariant of 3-manifolds.  The classical approach to remove the dependence
of~\eqref{eq:37} on the connection~$A$ is to treat~$S(A)$ as defining a
variational problem and to find its critical points.  That is indeed
interesting: the Euler-Lagrange equation asserts that the connection~$A$ is
flat.  There is not in general a unique flat connection, so no particular
critical value to choose as a topological invariant, though many interesting
topological invariants may be formed from the space of flat connections.
(They make an appearance in ~\S\ref{sec:7}.)  Witten's approach is quantum
mechanical \`a la Feynman: he integrates out the variable~$A$ to obtain a
topological invariant.  For each integer~$k$, termed the \emph{level} of the
theory, set
  \begin{equation}\label{eq:38}
     F_k(X) \text{``=''} \int_{\FX}e^{ikS(A)} \,dA. 
  \end{equation}
The integral takes place over a space of equivalence classes of connections.
This is the quotient~$\FX$ of the infinite-dimensional linear space~$\Omega
^1(X;\mathfrak{g})$ by a nonlinear action of the infinite-dimensional Lie
group of maps $X\to G$.  The action is given by the equation 
  \begin{equation*}
     g\cdot A = g\inv Ag + g\inv dg,\qquad g:X\to G,\quad A\in \Omega
     ^1(X;\mathfrak{g}). 
  \end{equation*}
The result of the integral is a complex number~$F_k(X)$.

I place quotation marks around the equality in~\eqref{eq:38} to indicate that
the path integral is only a heuristic: the presumptive measure~$dA$ in the
notation would have to be constructed before making sense of~\eqref{eq:38}.
Shortly I will indicate approaches to rigorously defining topological
invariants~$F_k(X)$.  Let me immediately point out a difficulty with the
naive formula~\eqref{eq:38}.  If there really were a measure~$dA$ which made
that formula work, then we would conclude that $F_k(X)$~is a topological
invariant of~$X$ which depends only on the orientation of~$X$.  After all,
the integrand---built from the classical Chern-Simons
invariant~\eqref{eq:37}---only depends on the orientation.  However, it turns
out that $F_k(X)$~depends on an additional topological structure on~$X$.  It
can variously be described as a 2-framing~\cite{A1}, a rigging~\cite{Se2}, or
a $p_1$-structure~\cite{BHMV}.  The dependence of the right hand side
of~\eqref{eq:38} on this structure signals the difficulties in defining the
integral.  Physicists approach the path integral through the process of
\emph{regularization}, which here can be accomplished by introducing a
Riemannian metric on~$X$.  The regularization does not preserve all of the
topological symmetry (orientation-preserving diffeomorphisms) of the
classical action.  One often says that the quantum theory has an
\emph{anomaly}.  Interestingly, to obtain a topological invariant, albeit of
a manifold with a $p_1$-structure, Witten introduces a ``counterterm'': the
Chern-Simons invariant of the auxiliary Riemannian metric, precisely the
invariant studied by Chern and Simons.  While anomalies are well-understood
geometrically in many other contexts, this Chern-Simons anomaly remains
somewhat of a mystery.
 
The implication of the previous may well be that one should construct a
measure~$dA$ on~$\field X$ which depends on a $p_1$-structure on~$X$.  As far
as I know, this has not been done.
 
There is an important extension of the invariant~\eqref{eq:38} to 3-manifolds
which contain an oriented link.  (Again there is an anomaly: the normal
bundle to the link must also be framed in the quantum theory.)  Suppose that
$\chi \:G\to\CC$ is a character of~$G$.  This means we represent~$G$ linearly
on a finite-dimensional vector space, and $\chi (g)$~is the trace of the
matrix which represents the action of~$g\in G$.  Given a connection~$A$
on~$X$ and an oriented loop~$C\subset X$ we define the \emph{holonomy} of~$A$
around~$C$ by solving a first-order ordinary differential equation.  The
holonomy is only defined up to conjugacy, but since the character is
invariant under conjugation we get a well-defined complex number~$\chi
_C(A)$.  Given a link $L=\amalg C_i$ which is the finite disjoint union of
oriented curves, the formal path integral definition is
  \begin{equation}\label{eq:41}
     F_k(X,L) \text{``=''} \int_{\FX}e^{ikS(A)} \prod\chi \mstrut _{C_i}(A)\,dA. 
  \end{equation}
(One may choose a different character~$\chi $ for each component of the
link.)  Witten shows how to compute the Jones polynomial invariant~\cite{Jo1}
of a link in~$S^3$ from the invariants~$F_k(X,L)$ for~$G=SU(2)$.  Indeed,
that was one of his motivations for this work.  The invariants
in~\eqref{eq:41} are a generalization of the Jones invariants to links in
arbitrary 3-manifolds.  They include as well the HOMFLYPT polynomial
invariants~\cite{FYHLMO}, \cite{PT}.
 
In another direction the Chern-Simons theory is defined for any compact Lie
group~$G$; the level~$k$ in~\eqref{eq:38} is replaced by a
class\footnote{There is a nondegeneracy restriction.  Namely $H^4(BG;\RR)$~
is naturally isomorphic to the vector space of $G$-invariant symmetric
bilinear forms on~$\Lie(G)$, and we require that the image of~$\lambda $ in
this vector space be a nondegenerate form.} $\lambda \in H^4(BG;\ZZ)$, where
$BG$~is the classifying space of the Lie group~$G$.  See the appendix for the
classical theory in this generality.  Note that $H^4(BSU(n);\ZZ)\cong \ZZ$
and there is a canonical generator.

How can we make mathematical sense out of~\eqref{eq:38}?  One approach might
be to attack directly the problem of constructing a measure~$dA$.  While this
has worked well in other quantum field theories, the theory at hand is in the
end topological and we might hope for a more topological approach.  So
instead we codify the structure inherent in~\eqref{eq:38} and seek to
construct examples of that structure.  In the remainder of this section we
elucidate the bare formal properties of the path integral.  This leads to a
set of axioms for a \emph{topological quantum field theory}, which we discuss
in~\S\ref{sec:1}.  In~\S\ref{sec:2} we return to Chern-Simons theories and to
their construction for various classes of Lie groups~$G$.
 
To begin, it is important to extend to compact 3-manifolds~$X$ with boundary.
In that case we replace~\eqref{eq:38} by
  \begin{equation}\label{eq:42}
     F_k(X)(\alpha ) \text{``=''} \int_{\FX(\alpha )}e^{ikS(A)} \,dA, 
  \end{equation}
where $\alpha $~is a connection on the boundary~$\bX$ and $\FX(\alpha )$~is
the space of gauge equivalence classes of connections on~$X$ whose boundary
value is~$\alpha $.  Thus $F_k(X)$~is a function on the space~$\FbX$ of gauge
equivalence classes on~$\bX$.   
 
\begin{figure} 
  \centering 
  \includegraphics[totalheight=120pt]{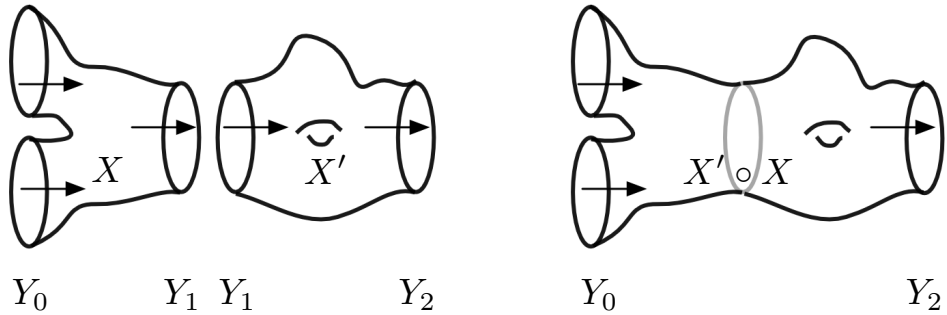} 
  \caption{Gluing of bordisms} \label{fig}
  \end{figure}

We abstract the formal properties of~\eqref{eq:42} for a general path
integral, not necessarily that in Chern-Simons theory.  A quantum field
theory lives in a fixed dimension~$n$, which is the dimension of the
``spacetimes'' in the theory.  For Chern-Simons theory $n=3$ and the theory
is topological, but in general we should imagine our manifolds as having
conformal structures or metrics, which we take to be Riemannian.  Physicists
call this ``Euclidean field theory''.  Certainly time is not present in this
positive definite signature---it is part of special relativity encoded in
Lorentz signature---so the use of `spacetime' in this Riemannian context
should be understood by analogy with the Lorentzian case.  Now for each
spacetime~$X$, which is an $n$-manifold with more structure, there is a
space~$\FX$ of fields and a classical action\footnote{As we see
from~\eqref{eq:38} it is only the exponentiated classical action which enters
into the path integral; indeed, in Chern-Simons theory only that exponential
is well-defined.  Also, in Euclidean field theory a minus sign usually
appears in place of~$i=\sqrt{-1}$ in~\eqref{eq:47} below, but we deliberately
use notation which follows the example of Chern-Simons.}
  \begin{equation*}
     S\:\FX\longrightarrow \CC. 
  \end{equation*}
The boundaries of spacetimes are closed $(n-1)$-manifolds~$Y$, and on these
we also imagine a space~$\FY$ of fields (but no classical action).  We
want~$Y$ to possess an infinitesimal normal bundle which is \emph{oriented}:
$Y$~comes equipped with an arrow of time.  In the Riemannian case there is
also a germ of a metric in the normal direction, but in a topological theory
the normal orientation is enough.  We consider a compact spacetime~$X$
together with a decomposition of its boundary as a disjoint union $\bX \cong
Y_0\amalg Y_1$.  The arrow of time on~$Y_0$ points into~$X$ and the arrow of
time on~$Y_1$ points out of~$X$.  Thus $Y_0$~is called the incoming part of
the boundary and $Y_1$~the outgoing part.  In topology we say $X$~is a
bordism from~$Y_0$ to~$Y_1$ and write $X\:Y_0\to Y_1$.  These bordisms
compose by gluing, which corresponds to the evolution of time.  (See
Figure~\ref{fig}.)

To a spacetime $X\:Y_0\to Y_1$ the ``semiclassical'' field theory attaches a
\emph{correspondence diagram}
  \begin{equation}\label{eq:44}
     \xymatrix@!C{&\field{X}\ar[dl]_s \ar[dr]^t \\ \field{Y_0} &&
     \field{Y_1}} 
  \end{equation}
The source map~$s$ and target map~$t$ are simply restriction of a field to
the appropriate piece of the boundary.  What is most important about the
fields is \emph{locality}. Suppose $X\:Y_0\to Y_1$ and $X':Y_1\to Y_2$ are
composable bordisms.  Then in the diagram
  \begin{equation}\label{eq:45}
     \xymatrix@!C{&&\field{X'\circ X}\ar[dl]_r \ar[dr]^{r'} \\
     &\field{X}\ar[dl]_s      \ar[dr]^t&&\field{X'}\ar[dl]_{s'} \ar[dr]^{t'}
     \\ \field{Y_0} && \field{Y_1} && 
     \field{Y_2}}   
  \end{equation}
the space of fields~$\field{X'\circ X}$ on the composition is the fiber
product of the maps~$t,s'$.  Loosely speaking, a field on the glued bordism
is a pair of fields on the separate bordisms which agree along the common
boundary.  One should keep in mind that fields are really infinite
dimensional \emph{stacks}.  In other words, fields may have automorphisms and
these play an important role in this context.  For example, in gauge theories
such as Chern-Simons theory the gauge transformations act as morphisms of
fields.  Certainly automorphisms must be accounted for in the maps and fiber
product in~\eqref{eq:44} and~\eqref{eq:45}.  The classical action is also
assumed local in the sense that
  \begin{equation}\label{eq:46}
  S_{X'\circ X}(\Phi )=S_X\bigl(r(\Phi ) \bigr) + S_{X'}\bigl(r'(\Phi ) \bigr) .
  \end{equation}
It is worth noting that whereas the locality of the classical action is
additive, that of the quantum theory is multiplicative.  Crudely,
quantization is a sort of exponentiation, as one can see from the integrand
in~\eqref{eq:42}.

Now we come to quantization in the geometric version given by the path
integral.  Assume there exist measures~$\mu\mstrut_X,\mu\mstrut_Y$ on the
spaces~$\field X,\field Y$.  Then define the Hilbert space~$\Hilb
Y=L^2(\field Y,\mu\mstrut_Y)$ and the linear map attached to a bordism
$X\:Y_0\to Y_1$ as the \emph{push-pull-with-kernel}~$e^{iS_X}$:
  \begin{equation}\label{eq:47}
     F_X=t_*\circ e^{iS_X}\circ s^*\:\Hilb{Y_0}\longrightarrow \Hilb{Y_1}. 
  \end{equation}
The pushforward~$t_*$ is integration.  Thus if $f\in
L^2(\field{Y_0},\mu\mstrut_{Y_0})$ and $g\in L^2(\field{Y_1},\mu\mstrut
_{Y_1})$, then
  \begin{equation*}
     \left\langle \bar{g},F_X(f) \right\rangle_{\Hilb{Y_1}} = \int_{\Phi \in
     \field X}\overline{g\bigl(t(\Phi ) \bigr)}\,f\bigl(s(\Phi )
     \bigr)\;e^{iS_X(\Phi )}\;d\mu\mstrut _X(\Phi ). 
  \end{equation*}
If the boundary of~$X$ is empty, then using the fact that $\field{\emptyset
}$~is a point, we see that \eqref{eq:47}~defines a complex number which is
the path integral~\eqref{eq:38}.
 
In short, the quantization~\eqref{eq:47} \emph{linearizes} the correspondence
diagram~\eqref{eq:44}.  This push-pull-with-kernel linearization of
correspondence diagrams is ubiquitous in mathematics.  The Fourier transform
is the archetypal example.  Let $S\:\RR\times \RR\to\RR$ be the duality
pairing $S(x,\xi )=x\xi $ and $s,t\:\RR\times \RR\to\RR$ the two projections.
Then the linearization of
  \begin{equation*}
     \xymatrix@!C{&{\RR\times \RR}\ar[dl]_s \ar[dr]^t \\ {\RR} &&
     {\RR}} 
  \end{equation*}
via~\eqref{eq:47} is the Fourier transform; the pushforward~$t_*$ is
integration with respect to suitably normalized Lebesgue measure.

The most important property of the path integral is \emph{multiplicativity}
under gluing.  It asserts that for a composition~$X'\circ X$ of bordisms as
in~\eqref{eq:45} we have
  \begin{equation*}
     F_{X'\circ X} = F_{X'}\circ F_X\:\Hilb{Y_0}\longrightarrow \Hilb{Y_2} 
  \end{equation*}
This is a formal consequence of the fiber product in~\eqref{eq:45} and the
additivity~\eqref{eq:46} if we assume an appropriate gluing law for the
measures.  As mentioned earlier, Chern-Simons can be defined for any compact
Lie group~$G$.  If $G$~is a finite group, then the pushforward~$t_*$
in~\eqref{eq:47} is a finite sum and the heuristics in this section can be
carried out rigorously~\cite{FQ}.  Also, for a finite group there are no
$p_1$-structures required and the quantum Chern-Simons theory is defined for
oriented manifolds. 
 
Finally, we remark that once quantum Chern-Simons theory is defined for
3-manifolds with boundary, then invariants of links are included.  Let $V$~be
the Hilbert space attached to the standard torus~$\cir\times \cir$.  In
Chern-Simons theory it is finite-dimensional, so we can fix a
basis~$e_1,e_2,\dots ,e_N$.  Suppose $X$~is an oriented 3-manifold which
contains an oriented link~$L=\amalg C_i$ and that the normal bundle to each
knot~$C_i$ is framed.  Let $X'$~be the compact 3-manifold with boundary
obtained from~$X$ by removing a tubular neighborhood of~$L$.  The orientation
and normal framing identifies each boundary component of~$X'$ with the
standard torus.  This identification is defined up to isotopy, and the
topological invariance of Chern-Simons theory implies that this is enough to
define an isomorphism of the Hilbert space of~$\partial X'$ with $\otimes _i
V$.  View~$\partial X'$ as incoming.  Then $F_{X'}\:\otimes _iV\to\CC$.
Suppose each component~$C_i$ of~$L$ is labeled by a basis vector~$e_i$.
Define
  \begin{equation*}
     F_k(X,L) = F_k(X')(\otimes _ie_i). 
  \end{equation*}
The constructions of Chern-Simons theory as a 1-2-3~theory produce a
canonical basis of~$V$ and assign a basis element to the characters~$\chi $
which appear in~\eqref{eq:41}.

  \section{axiomatization and two-dimensional theories}\label{sec:1}

Atiyah~\cite{A2} introduced a set of axioms for a topological quantum field
theory (TQFT) which is patterned after Segal's axiomatic treatment~\cite{Se2}
of conformal field theory.  Quinn~\cite{Q} gave a more elaborate treatment
with many elementary examples.  In recent years the concept of a TQFT has
broadened and the definition has evolved.  The most functorial modern
definition of a topological quantum field theory (TQFT) makes evident both
the analogy with homology theory and the formal properties of the path
integral outlined in~\S\ref{sec:5}.  Recall that homology is a functor
  \begin{equation*}
     H\:(Top,\amalg)\longrightarrow (Ab,\oplus ) 
  \end{equation*}
from the category of topological spaces and continuous maps to the category
of abelian groups and homomorphisms.  Both $Top$~and $Ab$~are \emph{symmetric
monoidal categories}: the tensor product on~$Top$ is disjoint union~$\amalg $
and on~$Ab$ is direct sum~$\oplus $.  The functor~$H$ is a symmetric monoidal
functor: homology is additive under disjoint union.  The codomain tensor
category may be replaced by, for example, the category~$Vect\mstrut _{\FF}$
of vector spaces over a field~$\FF$ where the monoidal structure is direct
sum of vector spaces.  We remark that a full list of axioms for homology
includes the Mayer-Vietoris property.
 
In a TQFT the category of topological spaces is replaced by a bordism
category of manifolds of fixed dimension.  As discussed in~\S\ref{sec:5} to
compose bordisms we need an ``arrow of time'', which we term a
\emph{collaring}.  Thus a collar of a closed manifold~$Y$ is an embedding
$Y\hookrightarrow \hY$ into a manifold~$\hY$ which is diffeomorphic
to~$(-\epsilon ,\epsilon )\times Y$ together with an orientation of the
normal bundle of~$Y\subset \hY$.  If $Y$~is a collared manifold, then the
opposite collared manifold~$-Y$ is the same embedding~$Y\subset \hY$ with the
reversed orientation on the normal bundle.  The boundary of a manifold~$X$
has a collar~\cite[\S4.6]{Hi}; the arrow of time points out of~$X$.  Let
$\bordO_n$ denote the category whose objects are compact collared
$(n-1)$-manifolds~$Y$.  A morphism $Y_0\to Y_1$ is a compact $n$-manifold~$X$
with a decomposition of its boundary $\bX = (\bX)_{\textnormal{in}}\;\amalg
\;(\bX)_{\textnormal{out}}$ and collar-preserving diffeomorphisms $-Y_0\to
(\bX)_{\textnormal{in}}$ and $Y_1\to(\bX)_{\textnormal{out}}$.
Bordisms~$X,X'$ define the same morphism if there exists a diffeomorphism
$X\to X'$ which commutes with the other data.  Composition is gluing of
bordisms.  Precise definitions are given in~\cite{GMWT}, for example.  The
disjoint union~$\amalg$ of manifolds endows~$\bordO_n$ with a monoidal
structure.

        \begin{definition}[]\label{thm:1}
 An \emph{$n$-dimensional topological quantum field theory}~$F$ is a
symmetric monoidal functor
  \begin{equation*}
     F\:(\bordO_n,\amalg )\longrightarrow (Ab,\otimes ). 
  \end{equation*}
        \end{definition}

\noindent
 Note that the tensor structure on abelian groups is \emph{tensor product}:
\emph{TQFT is multiplicative} whereas homology is additive.  As with homology
we can contemplate other codomains, for example replacing $(Ab,\otimes )$
with $(Vect_{\FF},\otimes )$ or $(R\textnormal{-}mod,\otimes )$ for some
commutative ring~$R$.  See Remark~\ref{thm:2} below for a slightly more
sophisticated replacement.  One can go all the way and simply declare that
the codomain is an arbitrary symmetric monoidal category, as in
Theorem~\ref{thm:29} below.  The domain $(\bordO_n,\amalg )$ may be replaced
by another bordism category: oriented bordism $(\bordSO,\amalg )$, spin
bordism $(\bordSpin,\amalg )$, framed bordism $(\bordfr,\amalg )$, etc.
 
The empty set is a manifold of any dimension; the empty manifold of
dimension~$n-1$ is a unit in~$(\bordO_n,\amalg )$.  It follows that
$F(\emptyset )=\ZZ$ for any TQFT~$F$.  A \emph{closed} oriented manifold~$X$
of dimension~$n$ is a morphism $X\:\emptyset \to\emptyset $ from which
$F(X)\in \ZZ$.  An $n$-dimensional TQFT therefore assigns to closed manifolds
algebraic objects as follows: 
  \begin{equation}\label{eq:4}
      \begin{tabular}{ c  c}  
      \toprule
      \emph{\vtop{\hbox{closed manifold}\hbox{\hskip 20pt of dim}}} &
      $F(\cdot )$\\  
      \midrule
      $n$ & element of $\ZZ$\\[-8pt]\\
      $n-1$ & \textnormal{$\ZZ$-module}\\ 
      \bottomrule
      \end{tabular}
  \end{equation}

We consider now \emph{oriented} TQFTs: the domain category is~$\bordSO_n$.
An object~$Y$ is a closed oriented $(n-1)$-manifold with a collar
$Y\hookrightarrow \hY$, and this induces an orientation on the
$n$-manifold~$\hY$.  The opposite manifold~$-Y$ has the opposite orientation
and the opposite collar.  The structure of a 1-dimensional oriented TQFT is
quite simple.  Let $F(\pt_{\pm})=A_{\pm}$ be the abelian groups assigned to a
point with orientation.  There is a unique closed oriented interval, up to
diffeomorphism, but it defines four distinct morphisms in~$\bordSO_1$:
  \begin{equation}\label{eq:5}
     \begin{aligned}  
      B = F(\,
     \vcenter{\hbox{\epsfysize=14pt\epsfbox{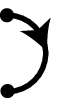}}  
     }\,) &\:A_+\otimes A_- \longrightarrow \ZZ \\
      B^{\vee}=F(\, 
     \vcenter{\hbox{\epsfysize=14pt\epsfbox{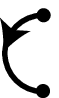}} }\,)
     &\:\ZZ\longrightarrow 
     A_-\otimes A_+ \\
      F(\,\vcenter{\hbox{\epsfysize=6pt\epsfbox{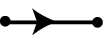}}
     }\,) &\:A_+\longrightarrow A_+ \\ 
      F(\,\vcenter{\hbox{\epsfysize=6pt\epsfbox{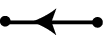}}
     }\,) &\:A_-\longrightarrow A_-\end{aligned} 
  \end{equation}
We read the pictures from left to right: the left boundary is incoming and
the right boundary outgoing.  It is easy to prove that the last two morphisms
are idempotents.  Assume they are identity maps; if necessary
replace~$A_{\pm}$ with the images of the idempotents.  Then an easy gluing
argument proves
  \begin{equation}\label{eq:6}
     \begin{alignedat}{2} A_+ &\xrightarrow{\id\times B^{\vee}} A_+\otimes
      A_-\otimes A_+ &&\xrightarrow{B\times \id} A_+ \\ A_- &\xrightarrow{
      B^{\vee}\times \id} A_-\otimes A_+\otimes A_- &&\xrightarrow{\id\times
      B} A_- 
      \end{alignedat}
  \end{equation}
are identity maps.  It follows that (i)\ $A_{\pm}$~are free and finitely
generated; (ii)\ $B$~determines an isomorphism $A_-\cong \Hom(A_+,\ZZ)$; and
(iii)\ $F(\cir)=\rank A_+$.  These arguments persist in an $n$-dimensional
theory: take the Cartesian product with a fixed $(n-1)$-manifold~$Y$.  

        \begin{proposition}[]\label{thm:10}
 Let $F$~be a TQFT defined on~$\bordSO_n$.  Then for any closed oriented
$(n-1)$-manifold~$Y$ the abelian group~$F(Y)$ is free and finitely generated,
there is a duality between~$F(Y)$ and~$F(-Y)$, and $F(\cir\times Y)=\rank
F(Y)$.

        \end{proposition}

        \begin{remark}[]\label{thm:2}
 Proposition~\ref{thm:10} already indicates the restrictive nature of
Definition~\ref{thm:1}: the abelian groups~$F(Y)$ are free.  To construct a
theory which includes torsion we can replace the target tensor category
$(Ab,\otimes )$ with the tensor category $(\textnormal{dg-}Ab,\otimes )$ of
differential graded abelian groups; we take the differential to have
degree~1.  The symmetry in $(\textnormal{dg-}Ab,\otimes )$ uses the Koszul
sign rule.  Then there is a 1-dimensional theory with
  \begin{equation}\label{eq:7}
     \begin{aligned} A_+\: &\cdots \longrightarrow 0\longrightarrow
       \ZZ\xrightarrow{\;n\;}\ZZ\longrightarrow 0\longrightarrow
     0\longrightarrow \cdots \\ 
     A_-\: 
       &\cdots \longrightarrow 0\longrightarrow0\longrightarrow
       \ZZ\xrightarrow{\;n\;}\ZZ\longrightarrow 0\longrightarrow \cdots
       \end{aligned} 
  \end{equation}
where the nonzero homogeneous groups in~$A_+$ are in degrees~$-1,0$ and
in~$A_-$ in degrees~$0,+1$.  Let $e_{-1},e_0$ and~$f_0,f_1$ be the obvious
basis elements in~\eqref{eq:7}. Set $B(e_0,f_0)=B(e_{-1},f_1)=1$ and
$B^{\vee}(1) = f_1\otimes e_{-1}+ f_0\otimes e_0$.  It is instructive to
check~\eqref{eq:6} and compute $F(\cir)=0$.  The complex~$A_+$ is a
resolution of~$\zn$: in this sense $F(\pt_+)$~is torsion.
        \end{remark}

There is a richer, but still simple, structure for~$n=2$.  Recall that a
\emph{Frobenius ring} is a unital ring~$A$ equipped with a homomorphism
$\theta \:A\to\ZZ$ which satisfies $\theta (xy)=\theta (yx)$ for all~$x,y\in
A$ and the pairing $x,y\mapsto \theta (xy)$ is nondegenerate.

        \begin{theorem}[]\label{thm:3}
 A 2-dimensional TQFT~$F$ determines a finitely generated commutative
Frobenius ring $A=F(\cir)$.  Conversely, a finitely generated commutative
Frobenius ring~$A$ determines a 2-dimensional TQFT~$F$ with $F(\cir)=A$. 
        \end{theorem}

\noindent
  The statement of Theorem~\ref{thm:3} can be found in many references, e.g.,
Dijkgraaf's thesis~\cite{D}.  A precise proof is given in~\cite{Ab}; see
also~\cite[\S A.1]{MS}.  A stronger statement---an equivalence of categories
of 2-dimensional TQFTs and finitely generated commutative Frobenius
rings---is proved in~\cite{K}.  We exhibit this theorem as the paradigmatic
\emph{generators-and-relations construction} of a TQFT.  The generators are the
unit, trace, multiplication, and comultiplication, which give 
  \begin{equation*}
     F(\quad \scalebox{.4}{\rput{90}(-.1,.25){\birthc}}\;),\qquad  F(\quad
     \scalebox{.4}{\rput{-90}(-.5,.25){\birthc}}\;),\qquad   F(\qquad
     \scalebox{.2}{\rput{90}(-.3,.5){\multc}}\;),\qquad F(\;\;
     \scalebox{.2}{\rput{-90}(-.3,.5){\multc}}\quad\;\;)
  \end{equation*}
respectively.  The commutativity and associativity relations assure that when
any bordism~$X$ is decomposed as a composition of these basic bordisms and
$F(X)$~is defined as a composition of the generating data, then $F(X)$~is
independent of the decomposition.

The notion of a TQFT, and the generators-and-relations Theorem~\ref{thm:3},
can be extended in various directions.  We can, for example, ask for
invariants of \emph{families} of manifolds in addition to single manifolds.
To see what this entails, suppose $F\:(\bordSO_n,\amalg )\to
(Vect_{\CC},\otimes )$ is an $n$-dimensional TQFT.  Let $\Y\to S$ be a fiber
bundle with fiber a closed oriented $(n-1)$-manifold.  Then the vector
spaces~$F(\Y_s)$ attached to fibers~$\Y_s$ fit together into a \emph{flat}
vector bundle we denote $F(\Y/S)\to S$.  Local trivializations and the flat
structure are derived from parallel transport: if $\gamma \:[0,1]\to S$ is a
smooth path from~$s_0$ to $s_1$, then $\gamma \inv \Y\:\Y_{s_0}\to \Y_{s_1}$
is a bordism and $F(\gamma \inv \Y)\:F(\Y_{s_0})\to F(\Y_{s_1})$ is defined
to be parallel transport.  Topological invariance shows that it is invariant
under homotopies of~$\gamma $.  If $\X\to S$ is a family of bordisms from
$\Y_0\to S$ to $\Y_1\to S$, then applying~$F$ we obtain
  \begin{equation*}
     F(\X/S)\in H^0\bigl(S;\Hom\bigl(F(\Y_0/S),F(\Y_1,S) \bigr) \bigr): 
  \end{equation*}
the topological invariance implies that $F(\X/S)$~is a flat section.  It is
natural, then, to extend~$F$ to a $\ZZ$-graded functor~$\Ftil$ which assigns
a cohomology class
  \begin{equation*}
     \Ftil(\X/S)\in H^\bullet\bigl(S;\Hom\bigl(F(\Y_0/S),F(\Y_1,S) \bigr) \bigr) 
  \end{equation*}
to a family $\X\to S$.  In addition to the functoriality under composition of
morphisms, we require naturality under base change.  This notion appears
for~$n=2$ in~\cite{KM} where the invariants are extended to nodal (better:
Deligne-Mumford stable) Riemann surfaces as well.  Gromov-Witten invariants
are an example of such a functor.  There is a large literature on this
subject: see~\cite{T} for a recent result about the structure of
2-dimensional theories in families.
 
In another direction we can ask that an $n$-dimensional TQFT extend to
manifolds of lower dimension, or higher codimension.  This idea dates from
the early '90s, when such theories were sometimes termed `extended
TQFTs'~\cite{L}, \cite{F1}.  They also go by the appellation `multi-tier
theories'~\cite{Se1}.  For $n=2$ a 3-tier theory is as far as we can go: it
assigns invariants to points.  We prefer to call this a `0-1-2 theory'.
Whatever the moniker, a 3-tier $n$-dimensional theory attaches a
\emph{category} to a manifold of codimension~2 (dimension~$n-2$).  If we work
over a ring~$R$, then the table~\eqref{eq:4} is extended to: 
  \begin{equation}\label{eq:10}
     \begin{tabular}{ c c} \toprule \emph{\vtop{\hbox{closed
      manifold}\hbox{\hskip 20pt of dim}}} & $F(\cdot )$\\ \midrule $n$ &
      element of $R$\\[-8pt]\\ $n-1$ & \textnormal{$R$-module}\\ [-8pt]\\
      $n-2$ & \textnormal{$R$-linear category}\\
      \bottomrule \end{tabular} 
  \end{equation}
(The morphism sets in an $R$-linear category are $R$-modules and composition
is an $R$-linear map.)  The entire structure is 2-categorical: the domain is
a 2-category whose objects are closed $(n-2)$-manifolds, 1-morphisms are
$(n-1)$-manifolds, and 2-morphisms are $n$-manifolds.  Families of manifolds
may also be incorporated using a topological version of a 2-category.  One
possibility is an \emph{infinity 2-category}, which is a simplicial set with
extra structure.  An infinity 1-category, also called a \emph{weak Kan
complex} or \emph{quasicategory}~\cite{BV}, \cite{Joy}, \cite{Lu}, may be
used to encode an ordinary TQFT for families of manifolds.  Other models,
such as Segal's $\Gamma $-spaces, may be more convenient.  The codomain can
be an arbitrary symmetric monoidal infinity 1- or infinity 2-category.

        \begin{remark}[]\label{thm:4}
 The discussion of ordinary 0-1 TQFTs surrounding~\eqref{eq:5} may be
restated in these terms: a TQFT $F\:(\bordSO,\amalg )\to \C$ with codomain a
symmetric monoidal category~$\C$ is determined up to isomorphism by a
\emph{dualizable object}~$A\in \C$.  Then $F(\pt_+)=A$ and the rest of the
theory is determined by specifying the duality data.  For~$\C=(Ab,\otimes )$
the dualizable objects are finitely generated free abelian groups.
        \end{remark}

        \begin{remark}[]\label{thm:33}
 One can also enrich homology theory to include categories analogously to
chart~\eqref{eq:10}.  But there is a very important difference with quantum
field theory.  The categories which occur in this way in homology theory are
groupoids---all morphisms are invertible.  In quantum field theories the more
general notion of duality replaces invertibility.
        \end{remark}

The structure of 0-1-2 theories, sometimes only partially defined, has
recently been elucidated by Costello~\cite{C}, Moore-Segal~\cite{MS}, and
others.  Building on this work Hopkins and Lurie~\cite{HL} prove a
generators-and-relations theorem for 0-1-2 theories, in essence giving
generators and relations for the infinity 2-category of 0-, 1-, and
2-manifolds.  (They also give an analogous result for 0-1~theories.)  The
simplest result to state is for \emph{framed} manifolds and was conjectured
by Baez-Dolan~\cite{BD}: the infinity~ 2-category of 0-, 1-, and 2-manifolds
with framing\footnote{Precisely, the framing is a trivialization of the
tangent bundle made 2-dimensional by adding a trivial bundle.  For example,
the framing of a point is a choice of basis of~$\RR^2$.} is freely generated
by a single generator, a framed point.  The precise statement is in terms of
functors to an arbitrary codomain.

        \begin{theorem}[Hopkins-Lurie~\cite{HL}]\label{thm:29}
 Let~$\C$ be a symmetric monoidal infinity 2-category.  Then the space of
0-1-2~theories of framed manifolds with values in~$\C$ is homotopy equivalent
to the space~$\C^{\textnormal{fd}}$ of \emph{fully dualizable} objects
in~$\C$.
        \end{theorem}

\noindent
 The notion of a fully dualizable object in a symmetric monoidal $n$-category
generalizes that of a dualizable object in a symmetric monoidal 1-category.
The statement for a 0-1-2 theory is simpler than for a 1-2 theory
(Theorem~\ref{thm:3}): there is now a single generator and no relation.
There is an explicit map which attaches to a field theory~$F$ the fully
dualizable object~$F(\pt_+)$, and this is proved to be a homotopy
equivalence.  One powerful feature of the theorem is its ``infinity~aspect'':
it applies to families so contains information about the diffeomorphism
groups.  The most important immediate application of this theorem is to
string topology~\cite{CS}; it implies that string topology produces homotopy
invariants of manifolds.  (Invariants from string topology which go beyond
the homotopy type are thought to exist if one extends to nodal
surfaces~\cite{Su}.)
 
There is a generalization of~Theorem~\ref{thm:29} to other bordism
categories.  The orthogonal group~$O_2$ acts on the infinity~2-category of
framed 0-, 1-, and~2-manifolds by rotating and reflecting the framing.
Therefore, by Theorem~\ref{thm:29} $O_2$~also acts on~$\C^{\textnormal{fd}}$.

        \begin{theorem}[]\label{thm:31}
 Let $G\to O_2$ be a homomorphism.  Then the space of 0-1-2~theories of
$G$-manifolds with values in~$\C$ is homotopy equivalent to the space of
homotopy $G$-fixed points in~$\C^{\textnormal{fd}}$.
        \end{theorem}

A 0-1-2 theory goes all the way down to points, so is as local as possible.
This explains why the structure of 0-1-2 theories is ultimately simpler than
that of 1-2 theories.  The principle that theories down to points are simpler
inspires the program of Stolz-Teichner~\cite{ST}, which aspires to construct
a quantum field theory model of the generalized cohomology theory known as
tmf~\cite{Ho}, which is closely related to elliptic cohomology theories.
These 3-tier theories also occur in two-dimensional \emph{conformal} field
theory: a category of ``D-branes'' is attached to a point in those theories.

Lurie has a detailed plan of proof for the conjectural generalization of
Theorem~\ref{thm:29} to higher dimensions.

        \begin{conjecture}[Baez-Dolan-Lurie]\label{thm:30}
  Let~$\C$ be a symmetric monoidal infinity $n$-category.  Then the space of
0-1-$\cdots$-$n$~theories of framed manifolds with values in~$\C$ is homotopy
equivalent to the space~$\C^{\textnormal{fd}}$ of \emph{fully dualizable}
objects in~$\C$.
        \end{conjecture}

\noindent
 The obvious generalization of Theorem~\ref{thm:31} is also conjectured to
hold.  The most local generators-and-relations constructions of 3-dimensional
Chern-Simons theory will rely on Conjecture~\ref{thm:30}.

  \section{three-dimensional tqft}\label{sec:2}

I do not know a structure theorem for 2-3 theories \`a la
Theorem~\ref{thm:3}.  Following the principle that a 3-tier theory is more
local than a 2-tier theory, and hence simpler, there \emph{is} a structure
theorem for 1-2-3 theories. It is proved for theories
defined\footnote{Work-in-progress by Davidovich, Hagge, and Wang aims to
prove that any 1-2-3~theory over~$\CC$ may be defined over a number field.}
over~$\CC$: the truncated 2-3~theory maps into~$(Vect_{\CC},\otimes )$.  An
early version of this result is contained in an influential manuscript of
Kevin Walker~\cite{Wa}, which builds on earlier work of Moore and
Seiberg~\cite{MSi}; there was also work of Kazhdan-Reshetikhin along these
lines and probably work of others as well.  A definitive version is in the
book of Turaev~\cite{Tu}, expanding on Reshetikhin-Turaev~\cite{RT}.

        \begin{theorem}[]\label{thm:5}
 A 1-2-3~theory~$F$ determines a \emph{modular tensor
category}~$\cat=F(\cir)$.  Conversely, a modular tensor category~$\cat$
determines a 1-2-3 theory~$F$ with $F(\cir)=\cat$.
        \end{theorem}

\noindent
 Loosely, a modular tensor category is the categorification of a commutative
Frobenius algebra.  More precisely, it is a braided monoidal category with
duals and a ribbon structure.  The category is required to be semisimple and
there is a nondegeneracy condition as well.  See~\cite{Tu} or~\cite{BK} for
details.  Theorem~\ref{thm:5} is again a generators-and-relations
construction: the monoidal structure tells $F(\qquad
\scalebox{.2}{\rput{90}(-.3,.5){\multc}}\;)$, the braiding tells $F$~applied
to the diffeomorphism of $\scalebox{.2}{\rput{90}(3,.5){\multc}}$ \quad\;
which exchanges the two incoming boundary components by a half-turn, the
duality tells $F$~applied to reflection on the circle, and the ribbon
structure is related to the generator of $\pi _1(\Diff^+\cir)\cong \pi
_1(SO_2)\cong \ZZ$.  The statement of Theorem~\ref{thm:5} is incomplete: we
have not told on which bordism category these theories are defined.  In fact,
the 1-, 2-, and 3-manifolds in the domain of~$F$ carry not only an
orientation but also a $p_1$-structure.  Together the orientation and
$p_1$-structure are almost a framing, but we do not
trivialize~$w_2$.\footnote{There are generalizations of Chern-Simons theory
defined on framed manifolds~\cite{J} and it will be interesting to have a
structure theorem analogous to Theorem~\ref{thm:5} for such 1-2-3~theories.}
 
The prime example and main object of interest for us is Chern-Simons theory,
which was introduced heuristically in~\S\ref{sec:5}.  Recall that the data
which defines the classical Chern-Simons invariant is a compact Lie group~$G$
and a class~$\lambda \in H^4(BG)$ termed the \emph{level}.  For $G$~a finite
group the path integral construction of a 2-3~theory is rigorous~\cite{DW},
\cite{FQ} as the integral in this case is a finite sum.  The path integral
was extended to higher codimensions in~\cite{F2} to construct a 1-2-3~theory
for $G$~finite.  In fact, the discussion in that work extends all the way
down to points, i.e., to a 0-1-2-3~theory, as we review presently.  These are
\emph{a priori} constructions of the Chern-Simons TQFT for finite gauge
groups.  For continuous gauge groups the only known constructions are with
generators and relations using Theorem~\ref{thm:5}.  The modular tensor
category~$\cat=\cat\gl$ should be the category of positive energy
representations of the loop group at level~$\lambda $ with the fusion
product, but I do not believe any detailed construction along these lines has
been carried out.  The main theorem in~\cite{FHT1} indicates that it should
also be a categorification of the twisted equivariant $K$-theory of~$G$ (with
$G$~action by conjugation and the twist derived from~$\lambda $), but this
too has not been carried out.  Rather, rigorous constructions---going back to
Reshetikhin and Turaev~\cite{RT}---for $G$~simply connected and simple take
as starting point the quantum group associated to the complexified Lie
algebra of~$G$ at a root of unity determined by~$\lambda $.  In particular,
this gives a rigorous construction of the $SU(n)$-quantum invariants~$F_k(X)$
introduced heuristically via the path integral in~\eqref{eq:42}.  There is a
separate line of development based on operator algebras; see the paper of
Jones~\cite{Jo2} in this volume for an account.  The theories for $G$~a torus
are quite interesting and have been classified by Belov-Moore~\cite{BM}.  The
corresponding modular tensor categories are special---every simple object is
invertible.  The relation between these modular tensor categories and toral
Chern-Simons theory is explored in~\cite{Sti}.  The classical toral theories
are specified by a finitely generated lattice with an even integer-valued
bilinear form, whereas the quantum theories only remember a quadratic form on
a finite abelian group extracted from the classical data.  Thus the map from
classical theories to quantum theories is many-to-one. 

        \begin{remark}[]\label{thm:34}
 A modular tensor category is semisimple.  As we saw in Remark~\ref{thm:2} we
need to alter the codomain to escape semisimplicity, even in a 2-tier theory.
The 3-dimensional Rozansky-Witten theory~\cite{RW}, \cite{Ro}, which takes as
starting data a complex symplectic manifold, is not semisimple so does not
fall under Theorem~\ref{thm:5}.  Ongoing work of Kapustin-Rozansky-Saulina
investigates the extension of Rozansky-Witten theory to a 0-1-2-3~theory.
        \end{remark}

A longstanding open question: Extend Chern-Simons theory down to points.  In
other words: Construct a 0-1-2-3~theory whose 1-2-3~truncation is
Chern-Simons theory for given~$\gl$.  A glance at the charts~\eqref{eq:4}
and~\eqref{eq:10} indicate the deep categorical waters ahead: a 4-tier
$n$-dimensional TQFT makes the assignments 
  \begin{equation*}
     \begin{tabular}{ c c c }  
      \toprule  
      \emph{\vtop{\hbox{closed manifold}\hbox{\hskip 20pt of dim}}}  
     & $F(\cdot )$  & category number\\  
      \midrule 
      $n$ & element of $R$ &$-1$\\[-8pt]\\  
      $n-1$ & $R$-module&$\hphantom{-}0$\\ [-8pt]\\
      $n-2$ & $R$-linear category&$\hphantom{-}1$\\  [-8pt]\\
      $n-3$ & $R$-linear 2-category&$\hphantom{-}2$\\  
      \bottomrule \end{tabular} 
  \end{equation*}
We attach a category number to each object in the chart: an $n$~category has
category number~$n$, a set has category number~$0$, and an element in a set
has category number~$-1$.\footnote{I once joked that every mathematician also
has a \emph{category number}, defined as the largest integer~$n$ such that
(s)he can think hard about $n$-categories for a half-hour without contracting
a migraine.  When I first said that my own category number was one, and in
the intervening years it has remained steadfastly constant whereas that of
many around me has climbed precipitously, if not suspiciously.}  It is a
feature of many parts of geometry over the past 25~years that the category
number of objects and theorems has increased.  Whereas theorems about
equivalence classes---sets---used to be sufficient, new questions demand that
automorphisms be accounted for: whence categories.  This trend has
affected---some would say infected---parts of quantum field theory as well.
In our current context it appears that a 0-1-2-3~TQFT attaches a dualizable
object in a symmetric monoidal 3-category to an oriented point.

Fortunately, it may be possible to extend Chern-Simons down to points using
only ordinary 1-categories.  As motivation, one categorical level down we
observe that a ring~$R$ determines a 1-category, the category of $R$-modules.
Similarly, a \emph{monoidal} 1-category~$\R$ determines a 2-category of its
modules.  We can hope that in a given 0-1-2-3~theory~$F$, such as
Chern-Simons theory, the 2-category ~$F(\pt_+)$ is the category of
$\R$-modules for a monoidal 1-category~$\R$.  We further speculate that if
so, then $F(\cir)$~is the \emph{Drinfeld center}~$Z(\R)$
of~$\R$.\footnote{One can prove with reasonable hypotheses that
$F(\cir)$~\emph{is} the Hochschild homology of~$\R$; the center is the
Hochschild cohomology, and in favorable cases these may be identified.}  The
Drinfeld center of any monoidal category is braided, and in favorable
circumstances~\cite{Mu} it is a modular tensor category.  As evidence in
favor of these hypotheses in the case of Chern-Simons theory, we
exhibit~$\R$ in case the gauge group~$G$ is finite  and show
that the modular tensor category~$F_{\gl}(\cir)$ is its Drinfeld center.

To begin we recall the precise definition of the center.
        \begin{definition}\label{thm:6}
 Let $\R$~be a monoidal category.  Its \emph{Drinfeld center}~$Z(\R)$ is the
category whose objects are pairs~$(X,\epsilon )$ consisting of an object~$X$
in~$\R$ and a natural transformation $\epsilon (-)\:X\otimes -\to -\otimes
X$.  The transformation~$\epsilon $ is compatible with the monoidal structure
in that for all objects~$Y,Z$ in~$\R$ we require 
  \begin{equation*}
     \epsilon (Y\otimes Z)= \bigl(\id_Y\otimes \epsilon (Z) \bigr)\circ
     \bigl(\epsilon (Y)\otimes \id_Z \bigr). 
  \end{equation*}
        \end{definition} 
 \noindent This notion of center was introduced by Joyal-Street~\cite{JS}.  A
recent paper by Ben-Zvi, Francis, and Nadler~\cite{BFN} extends this notion
of a center to a much broader context and contains discussions pertinent to
this paper.

We explicate~(ii).  Suppose $G$~is a finite group and $\lambda \in H^4(BG)$ a
level.  There is an \emph{a priori} ``path integral'' construction of the
associated quantum Chern-Simons theory~$F_{\gl}$ in~\cite{F2} which starts
with the classical Chern-Simons invariant for a finite group.  The fields in
this theory are principal $G$-bundles which, since $G$~is finite, are
covering spaces with Galois group~$G$.  Over the circle the groupoid of
principal $G$-bundles is equivalent to the groupoid~$G\gpd G$ of $G$~acting
on itself by conjugation.  The level~$\lambda $, through classical
Chern-Simons theory, produces a central extension of this groupoid: for every
pair of elements~$x,y\in G$ a hermitian line~$L_{x,y}$, for every triple of
elements~$x,y,z\in G$ an isomorphism
  \begin{equation}\label{eq:13}
     L_{yxy\inv ,z}\otimes L_{x,y}\longrightarrow L_{x,zy} 
  \end{equation}
and a consistency condition on the isomorphisms~\eqref{eq:13} for quartets of
elements.  Then $F_{\gl}(\cir)$~is the category of $L$-twisted
$G$-equivariant vector bundles over~$G$.  We remark that this is a concrete
model for \emph{twisted $K$-theory}, where the twisting is defined by~$L$.
(The relation of Chern-Simons theory to twisted $K$-theory is discussed
in~\S\ref{sec:6}.)  The monoidal structure is convolution, or pushforward
under multiplication.  We now exhibit this monoidal category as the center of
a monoidal tensor category~$\R_{\gl}$.
 
First, view~$H^4(BG)$ as $H^2(BG;\CP^{\infty})$ and so~$\lambda $ as
representing a central extension of~$G$ by the abelian group-like (Picard)
category of hermitian lines.  A cocycle\footnote{The construction, and so the
theory, depends on the choice of cocycle up to noncanonical isomorphism.  In
the notation `$\lambda $'~should be understood to include the choice of
cocycle.}  which represents~$\lambda $ is then a hermitian line~$K_{x,y}$ for
every pair of elements~$x,y\in G$, cocycle isomorphisms
  \begin{equation*}
     K_{x,y}\otimes K_{xy,z}\longrightarrow K_{x,yz}\otimes K_{y,z} 
  \end{equation*}
for triples of elements~$x,y,z\in G$, and a consistency condition on these
isomorphisms for quartets of elements in~$G$.  This is a categorified version
of the usual cocycle for a central extension of a discrete group by the
circle group.  We normalize $K_{1,x}=K_{x,1}=\CC$ for all~$x\in G$.  Let
$\R_{\gl}$~be the category of complex vector bundles on~$G$, or equivalently
the category of $G$-graded complex vector spaces.  If $W=\{W_x\}_{x\in G}$
and $W'=\{W'_y\}_{y\in G}$ are objects of~$\R_{\gl}$, set
  \begin{equation*}
     (W\otimes W')_z = \bigoplus\limits_{xy=z} K_{x,y}\otimes W_x\otimes
     W'_y. 
  \end{equation*}
This monoidal structure on~$\R_{\gl}$ is a twisted convolution, and so
$\R_{\gl}$~is a twisted group algebra of~$G$ with coefficients in the
category of complex vector spaces.  The dual~$W^*$ of an object~$W$ is
defined as 
  \begin{equation*}
     (W^*)_x = K^*_{x,x\inv }\otimes (W_{x\inv })^*. 
  \end{equation*}
Then it is straightforward to check from Definition~\ref{thm:6} that an
object in the center~$Z(\R_{\gl})$ of~$\R_{\gl}$ is a vector bundle $W\to G$
together with isomorphisms 
  \begin{equation}\label{eq:17}
     L_{x,y}\otimes W_x\longrightarrow W_{yxy\inv } \;,
  \end{equation}
where 
  \begin{equation}\label{eq:18}
     L_{x,y} := K^*_{yxy\inv ,y}\otimes K_{y,x}\;. 
  \end{equation}
Then \eqref{eq:13} follows from~\eqref{eq:18} and \eqref{eq:17}~satisfies a
consistency condition compatible with~\eqref{eq:13}.  This gives rise to a
functor 
  \begin{equation}\label{eq:19}
     Z(\R_{\gl})\longrightarrow F_{\gl}(\cir)\;, 
  \end{equation}
and it is not hard to check that \eqref{eq:19}~is an equivalence of braided
monoidal categories, in fact of modular tensor categories.

        \begin{remark}[]\label{thm:8}
 This construction appears in~\cite{F2}---see especially~(9.16)---though not
in this language.  The arguments in that reference make clear that
$F_{\gl}$~extends to a 0-1-2-3~theory with $F_{\gl}(\pt_+)=\R_{\gl}\module$. 
        \end{remark}

        \begin{remark}[]\label{thm:9}
 The modular tensor category~$F_{\gl}(\cir)$ is the category of
representations of a Hopf algebra which is the Drinfeld double of a group
algebra~$\CC[G,\lambda ]$.  The equivalence~\eqref{eq:19} is a categorified
version: $F_{\gl}(\cir)$~is the Drinfeld center of the twisted group
algebra~$\R_{\gl}$ of~$G$ with coefficients in~$Vect_{\CC}$.
        \end{remark}

        \begin{remark}[]\label{thm:28}
  A variant of this construction works when $G$~is a torus group~\cite{FHLT}.
Also, a different approach to Chern-Simons for a point using conformal nets
is being developed by Bartels-Douglas-Henriques.
        \end{remark}

  \section{dimensional reduction and integrality}\label{sec:6}

A quantum field theory in $n$-dimensions gives rise to theories in lower
dimension through a process known as \emph{dimensional reduction}.  In
\emph{classical} field theory as practiced by physicists, fields in an
$n$-dimensional theory are functions on $n$-dimensional Minkowski
spacetime~$M^n$.  The entire theory is invariant under the Poincar\'e
group~$P^n$, which acts on~$M^n$ by affine transformations.  Dimensional
reduction in this context~\cite[\S2.11]{DF} is carried out by restricting the
theory to fields which are invariant under a $k$-dimensional vector space~$T$
of translations.  Such fields drop to the quotient affine space~$M^n/T$,
which may be identified as~$M^{n-k}$.  The dimensionally reduced theory is
invariant under the Poincar\'e group~$P^{n-k}$.  But there is an additional
symmetry: the rotation group~$SO(k)$ of the space~$T$.  In other words,
\emph{the dimensional reduction of a theory has extra symmetry}.  In this
section we explain a conjectural analog for topological quantum field theories.
 
Dimensional reduction may be carried out in \emph{quantum} field theory as
well.  In the axiomatic framework this is easy, and for convenience we
restrict to the dimensional reduction of a 3-dimensional TQFT~$F$ to a
2-dimensional TQFT~$F'$.  The definition is quite simple: for any
manifold~$M$ set
  \begin{equation*}
     F'(M) = F(\cir\times M). 
  \end{equation*}
It is easy to check that this defines a 2-dimensional theory~$F'$.  Suppose
$F$~is a theory defined over the complex numbers; that is true for
Chern-Simons theory, which is our main example. 

        \begin{observation*}[]
 For a closed 2-manifold~$Y$ the complex number~$F'(Y)$ is an integer.  
        \end{observation*}

\noindent
 This follows easily as $F'(Y)=F(S^1\times Y)=\dim F(Y)$ is the dimension of
a complex vector space, hence is a nonnegative integer.\footnote{We can
consider theories in which the complex vector spaces are $\zt$-graded.  In
that case we would not have nonnegativity of~$F'(Y)$.}  We record this in a
chart:
  \begin{equation}\label{eq:52}
     \begin{tabular}{ c c c} \toprule \emph{\vtop{\hbox{closed
      manifold}\hbox{\hskip 20pt of dim}}} & $F(\cdot )$& $F'(\cdot )$\\
     \midrule $3$ &element of $\CC$\\[-8pt]\\ $2$ & \textnormal{$\CC$-vector
     space} & \textnormal{element of~$\ZZ$}\\ [-8pt]\\ 
      $1$ & \textnormal{$\CC$-linear category}& \textnormal{$\ZZ$-module}\\
     \bottomrule \end{tabular}  
  \end{equation}
The last entry in the lower right hand corner is speculative, as we do not
know that the complex vector space $F(S^1\times S^1)$ can be refined to an
abelian group.  Nonetheless, we make the following

        \begin{conjecture}[]\label{thm:16}
 The dimensional reduction~$F'$ is defined over~$\ZZ$.   
        \end{conjecture}

\noindent
 In other words, there is a functor from a bordism category to~$(Ab,\otimes
)$ whose composition with the functor~$-\otimes \mstrut _{\ZZ}\CC\:(Ab,\otimes
)\to (Vect\mstrut _{\CC},\otimes )$ is~$F'$.  

To formulate Conjecture~\ref{thm:16} precisely requires that we specify the
precise nature of the linear category attached to a 1-manifold by the
1-2-3~theory~$F$.  We will not attempt that here, but report that with
appropriate ``smallness'' hypotheses one can show that $F(S^1\times S^1)$~is
the \emph{Hochschild homology} of the linear category~$F(S^1)$.  Therefore,
we seek an abelian group~$A$ which refines the Hochschild homology in the
sense that the latter is naturally~$A\otimes\mstrut _\ZZ \CC$.  This
immediately brings to mind \emph{$K$-theory}.  Fortuitously, a similar
integrality is required in noncommutative Hodge theory~\cite[\S2.2.6]{KKP}.
There is a well-formulated plan, due to Bondal and To\"en, to define the
appropriate $K$-theory of appropriate dg-categories, based on work of
To\"en-Vaquie~\cite{TV} and the semi-topological $K$-theory of
Friedlander-Walker~\cite{FW}.\footnote{In fact, this $K$-theory group would
refine the periodic cyclic homology, not the Hochschild homology, but these
are closely related.}

Now that $K$-theory is on the table we may be bolder and ask that $F'$~map
into the world of stable homotopy theory, which is a natural home for
$K$-theory.  In that world the integers are replaced by the sphere
\emph{spectrum}~$S$, so at first we can replace~`$\ZZ$' in the last column
of~\eqref{eq:52} with~`$S$'.  Modules over the sphere spectrum are
spectra---just as modules over the integers are abelian groups---so we would
have a TQFT with values in spectra.  A more refined conjecture is that we can
define~$F'$ to have values in $K$-modules, where `$K$'~denotes $K$-theory, an
$E^{\infty}$ ring spectrum.  In other words, we refine
Conjecture~\ref{thm:16} to

        \begin{conjecture}[]\label{thm:17}
 The dimensional reduction~$F'$ is defined over~$K$. 
        \end{conjecture}

\noindent
 We remark that spectra have appeared in the context of Floer
homology~\cite{CJS}, \cite{M}, and most likely this fits into TQFT and
Conjecture~\ref{thm:17}.  Also, this refined conjecture is made in the
noncommutative Hodge theory.  Finally, one would naturally extend these
speculations to propose a 0-1-2 theory~$F'$ over~$K$.  This seems quite
mysterious, however, even in the explicit example of Chern-Simons for a
finite group.
 
While these speculations are rather abstract, there is a theorem which
provides some support.  Looked at the other way, these speculations provide a
context for the theorem, and indeed closely related considerations motivated
it in the first place~\cite{F4}.  Namely, consider the case of Chern-Simons
theory~$F$ for a compact Lie group~$G$ and level $\lambda \in H^4(BG;\ZZ)$.
Denote the reduction to a 2-dimensional theory by~$F'$.  Then
Conjecture~\ref{thm:17} suggests that $F'(\cir)$~is a $K$-theory group and
can even be refined to a spectrum which is a $K$-module.  Sticking for the
moment to abelian groups, notice from Theorem~\ref{thm:3} that $F'(\cir)$ is
a Frobenius ring.  In this context physicists call this the \emph{Verlinde
ring}, which first appeared in rational conformal field theory~\cite{V}.  Let
us denote it as~$R(G,\lambda )$.

        \begin{theorem}[F.-Hopkins-Teleman]\label{thm:18}
 There is a natural isomorphism 
  \begin{equation}\label{eq:53}
     \Phi \:R(G,\lambda )\longrightarrow K_G^{\tau (\lambda )}(G). 
  \end{equation}
        \end{theorem}

\noindent
 The right hand side is the twisted equivariant $K$-theory of the Lie
group~$G$, where $G$~acts on itself by conjugation.  The twisting~$\tau
(\lambda )$ is transgressed from~$\lambda $ and then shifted by a constant
twisting.\footnote{This is the well-known shift $k\to k+n$ in the
$SU(n)$-theories.}  See~\cite{FHT1} for a detailed development and proof.

        \begin{remark}[]\label{thm:19}
 In Theorem~\ref{thm:18} the left hand side $R(G,\lambda )$~is defined to be
the free abelian group generated by positive energy representations of the
loop group of~$G$ at a fixed level; the product is the fusion product.  This
fits the Chern-Simons story if we construct the modular tensor
category~$F(\cir)$ as the category of these representations.  The theorem
itself lies at the intersection of representation theory and algebraic
topology.
        \end{remark}

        \begin{remark}[]\label{thm:20}
 While the statement of Theorem~\ref{thm:18} was motivated by this physics,
the explicit map~$\Phi $ was not.  (There is other, closely related, physics
which could serve as motivation: the supersymmetric WZW model~\cite{Mi}.)
The explicit formula for~$\Phi $ justifies its name: the \emph{Dirac family}
associated to a representation of the loop group.  There is an analog for
representations of the finite dimensional group~$G$.  These Dirac families
should find more applications in representation theory.
        \end{remark}

        \begin{remark}[]\label{thm:21}
 There is a refinement over~$K$ in this case: one can make a $K$-module whose
nontrivial homotopy group is $K_G^{\tau (\lambda )}(G)$. 
        \end{remark}

        \begin{remark}[]\label{thm:22}
 While the 3-dimensional Chern-Simons theory~$F$ is defined on a bordism
category of oriented manifolds with $p_1$-structure, the 2-dimensional
reduction~$F'$ factors to a theory of oriented manifolds: no
$p_1$-structures.
        \end{remark}

The known rigorous constructions of the 3-dimensional Chern-Simons theory~$F$
use generators and relations---for the 1-2-3 theory as encoded in a modular
tensor category (Theorem~\ref{thm:5}) and for the 0-1-2-3 theory at least
conjecturally in some cases (Conjecture~\ref{thm:30}, Remark~\ref{thm:28}).
One attractive feature of the 2-dimensional reduction~$F'$ is that it has an
\emph{a priori} construction~\cite{FHT2} which follows the path integral
heuristic discussed in~\S\ref{sec:5}.  Recall that if $Y\:S_0\to S_1$ is a
2-dimensional bordism, then the path integral linearizes the correspondence
diagram
  \begin{equation*}
     \vcenter{\xymatrix@!C{&\field{Y}\ar[dl]_s \ar[dr]^t \\ \field{S_0} &&
     \field{S_1}} }\qquad \qquad  t_*\circ e^{iS_Y}\circ
     s^*\:L^2(\field{S_0})\longrightarrow  L^2(\field{S_1}) 
  \end{equation*}
Here $\mathcal{F}$~is the \emph{infinite-dimensional} stack of
$G$-connections.  The definition of the $L^2$-spaces and the
pushforward~$t_*$ relies on measures which are consistent under gluing (and
which do not usually exist, which is why this is a heuristic).  Our \emph{a
priori} construction replaces this path integral by a topological version: 
  \begin{equation*}
     \vcenter{\xymatrix@!C{&\fconn{Y}\ar[dl]_s \ar[dr]^t \\ \fconn{S_0} &&
     \fconn{S_1}} }\qquad \qquad t_*\circ  s^*\:K(\fconn{S_0})\longrightarrow
     K(\fconn{S_1})  
  \end{equation*}
In this topological diagram the infinite-dimensional stack~$\mathcal{F}$ of
$G$-connections has been replaced by the \emph{finite-dimensional}
stack~$\mathcal{M}$ of \emph{flat} $G$-connections.  The classical action~$S$
vanishes.  Also, the Hilbert spaces of $L^2$-functions are replaced by the
abelian $K$-theory groups.  The push-pull linearization is pure topology: the
pullback~$s^*$ is functorially defined and the pushforward~$t_*$ depends on
consistent $K$-theory orientations of the maps~$t$.  The consistency
conditions refer to gluing of bordisms as expressed in diagram~\eqref{eq:45},
but with the moduli stacks~$\mathcal{M}$ of flat $G$-connections replacing
the stacks~$\mathcal{F}$ of all $G$-connections.  It is important that flat
$G$-connections are local in the sense that $\fconn{X'\circ X}$ is a fiber
product, just as $\field{X'\circ X}$~is.  This topological linearization of a
correspondence diagram brings us full circle back to classical topology,
though here we use $K$-theory in place of ordinary cohomology.
 
The moduli stack~$\fconn{S^1}$ of (flat) $G$-connections on the circle is
equivalent to the global quotient~$G\gpd G$ of $G$~acting on itself by
conjugation.  Therefore, $K(\fconn{\cir})\cong K(G\gpd G)\cong K_G(G)$.  A
consistent orientation includes a twisting of this $K$-theory group, and in
fact it is the group which occurs on the right hand side of~\eqref{eq:53}.
This construction induces a Frobenius ring structure on this twisted
equivariant $K$-theory group: the product is analogous to the Pontrjagin
product on the homology of a topological group.

        \begin{remark}[]\label{thm:23}
 In~\cite{FHT2} we prove that consistent orientations exist and are induced
from a certain ``universal orientation''.  The group of universal
orientations maps to the group of levels, which is~$H^4(BG;\ZZ)$, but in
general that map is neither injective nor surjective.  (It is a bijection for
connected and simply connected groups such as~$SU(n)$.)  As the
2-dimensional reduction of Chern-Simons theory depends on a universal
orientation, we can ask whether that is true of the whole three-dimensional
Chern-Simons theory, which is usually thought to depend only on the level. 
        \end{remark}

        \begin{remark}[]\label{thm:24}
 We emphasize the analogy between measures and orientations.  In our formal
picture of the path integral measures enable analytic integration.  In this
topological construction orientations enable topological integration.
        \end{remark}

  \section{summary and Broader perspectives}\label{sec:7}

The current state of affairs for constructions of quantum Chern-Simons
theories is:

 \begin{itemize}
 \item There is a generators-and-relations construction of the 1-2-3 theory
via modular tensor categories for many classes of compact Lie groups~$G$.
These include finite groups, tori, and simply connected groups, the latter
via quantum groups or operator algebras. (\S\ref{sec:2})
 \item There are new generators-and-relations constructions---at this stage
still conjectural---of the 0-1-2-3 theory for certain groups, including
finite groups and tori. (\S\ref{sec:2})
 \item There is an \emph{a priori} construction of the 0-1-2-3 theory for a
finite group. (\S\ref{sec:2})
 \item There is an \emph{a priori} construction of the dimensionally reduced
1-2 theory for \emph{all} compact Lie groups~$G$. (\S\ref{sec:6})
 \end{itemize}

\noindent
 In many ways this represents a strong understanding of this particular
three-dimensional topological field theory.  It is a small part of the
successes mathematicians have achieved in understanding QFT-Strings over the
past 25~years.  On the credit side too are the general structure theorems for
topological field theories recounted in~\S\ref{sec:1} and~\S\ref{sec:2}.  All
of this may be viewed as a triumph of the axiomatization of path integrals,
which was motivated in~\S\ref{sec:5} and carried out in the succinct
Definition~\ref{thm:1}.  But I want to argue now that, in fact, this
axiomatization falls far short of capturing all properties of the path
integral.  In particular, it fails for the Chern-Simons path integral, which
is purely topological even at the classical level.

        \begin{remark}[]\label{thm:27}
 The beauty of the path integral heuristic is that it is an \emph{a priori}
definition of the topological invariants.  Indeed, one of Witten's
motivations in~\cite{W1} was to find a description of the Jones polynomial
invariant of a link which manifestly exhibits 3-dimensional topological
invariance.  The bullet points above do not include an \emph{a priori}
rigorous construction of these invariants except for a finite gauge group, in
which case the path integral heuristic can be expanded and made rigorous.
 
Another observation: the path integral takes as input the classical
Chern-Simons invariant.  But none of the mathematical constructions described
in the bulleted list above uses this invariant save the \emph{a priori}
construction for a finite group.
        \end{remark}

To illustrate our lack of understanding let us return to the Chern-Simons
path integral~\eqref{eq:38} on a closed 3-manifold~$X$:
  \begin{equation}\label{eq:56}
     F_k(X) \text{``=''} \int_{\FX}e^{ikS(A)} \,dA. 
  \end{equation}
Now the usual expression for a path integral has integrand $e^{iS/\hbar}$,
where $\hbar$~is Planck's constant.  It is a \emph{constant} of nature:
$\hbar \sim 1.054572 \times 10^{-34}\; \text{m}^2\, \text{kg} / \text{sec}$.
But in this context we can consider the limit~$\hbar\to0$ in which quantum
effects are suppressed, the so-called semiclassical limit.  Mathematically,
an oscillatory integral of the form
  \begin{equation}\label{eq:57}
     F(\hbar) = \int e^{iS(x)/\hbar}\;dx 
  \end{equation}
has an asymptotic expansion as~$\hbar\to0$ which is derived by the method of
stationary phase~\cite{B}.  It is a sum over the critical points of~$S$ if
$S$~is a Morse function or an integral over the critical manifold if $S$~is
Morse-Bott; in physics this sum is expressed in terms of Feynman diagrams and
is the basis of perturbation theory and the many spectacular computations in
quantum field theory.  Comparing~\eqref{eq:56} and~\eqref{eq:57} we see that
the level~$k$ plays the role of~$1/\hbar$.  In other words, in this
topological theory $\hbar$~takes on a countable set of values tending to
zero.  So the semiclassical limit is~$k\to\infty $ and we expect an
asymptotic expansion which is a sum over the critical points of~$S$.  Recall
that the critical points of~$S$ are the flat connections on~$X$.
Choose~$G=SU(2)$ and suppose that $X$~is a closed oriented 3-manifold on
which there is a finite set of equivalence classes~$\fconn{X}$ of flat
$SU(2)$-connections.  Let $\fconn{X}^0$~denote the equivalence classes of
\emph{irreducible} flat $SU(2)$-connections, assumed nonempty.  These
contribute to the leading order asymptotics of~$F_k(X)$ as $k\to\infty $
(see~\cite{W1}, \cite{FG}).

        \begin{conjecture}[]\label{thm:32}
 The quantum Chern-Simons invariant satisfies
  \begin{equation}\label{eq:58}
     F_k(X)\sim \frac 12\;e^{-3\pi i/4}\sum\limits_{A\in \fconn{X}^0}
     e^{2\pi iS_X(A)(k+2)}\,e^{-2\pi i I_{A}/4}\,\sqrt{\tau _X(A)}\,. 
  \end{equation}
        \end{conjecture}

\noindent
 In~\eqref{eq:58} $S_X(A)$~is the classical Chern-Simons
invariant~\eqref{eq:37}, $I_X(A)$~is the Atiyah-Patodi-Singer spectral
flow~\cite{APS}, and $\tau _X(A)$~is the Franz-Reidemeister torsion~\cite{R},
\cite{Fra}.
 
Now to the point: The left hand side is rigorously defined by the explicit
construction of Chern-Simons theory using quantum groups.  The right hand
side is a sum of classical invariants of flat connections.  Thus
\eqref{eq:58}~is a well-formulated mathematical statement, a conjecture
derived from the path integral.  I offer the fact that this most basic
consequence of the path integral is not proved as evidence that the
axiomatics of~\S\ref{sec:1} do not capture all of its essential features.  I
will return to this line of thought presently.

        \begin{remark}[]\label{thm:25}
 The left hand side of~\eqref{eq:58} is defined from the theory of quantum
groups or loop groups.  The computation of the exact quantum invariant, say
for Seifert fibered manifolds, has as a crucial ingredient explicit formulas
in that literature~\cite{KaWa}.  On the other hand, the right hand side
involves three invariants of flat connections which are seemingly unrelated
to quantum groups or loop groups.  Conjecture~\ref{thm:32} is typical of many
from quantum field theory and string theory in that it relates parts of
mathematics hitherto unconnected.
        \end{remark}

        \begin{remark}[]\label{thm:26}
 There are some special cases of Conjecture~\ref{thm:32} which have been
proved in the mathematics literature.  The proofs go by examining the
explicit formula for~$F_k(X)$ and using techniques of classical analysis to
derive the large~$k$ asymptotics.  For certain lens spaces this was carried
out by Jeffrey~\cite{Jef} and Garoufalidis~\cite{Gar}.  Rozansky~\cite{Roz}
generalized to all Seifert fibered manifolds.  Of course, most 3-manifolds
are not Seifert fibered, and for these \eqref{eq:58}~remains a conjecture.
        \end{remark}

  \begin{table}[h]
     \begin{tabular}{ c@{\qquad } *{3}{r@{\,}c@{\,}l@{\qquad}} } 
     \toprule $k$ & \multicolumn{3}{c}{exact value $F_k(X)$\;\;\;}&  
       \multicolumn{3}{c}{asymptotic value\;\;\;}&  
       \multicolumn{3}{c}{ratio\;\;\;\;\;\;}\\
     \midrule 
       141 & 0.607899& $+$& 0.102594\,i& 0.596099& $+$& 0.151172\,i& 
         0.999182& $-$& 0.081285\,i\\ [-8pt] \\
       142 & $-$0.104966& $-$& 0.151106\,i& $-$0.094614& $-$& 0.157913\,i
         & 0.997181& $-$& 0.067244\,i\\ [-8pt] \\
       143 & 0.123614& $-$& 0.139016\,i& 0.132261& $-$& 0.128045\,i& 
         1.007707& $-$& 0.075491\,i\\ [-8pt] \\
       144 & $-$0.612014& $+$& 0.038199\,i& $-$0.614913& $-$& 0.008261\,i
         & 0.994271& $-$& 0.075479\,i\\ [-8pt] \\
       145 & $-$0.291162& $-$& 0.132171\,i& $-$0.281928& $-$& 0.153204\,i& 
         0.993986& $-$& 0.071336\,i\\ [-8pt] \\
       146 & $-$0.413944& $+$& 0.674785\,i& $-$0.465909& $+$& 0.642185\,i& 
         0.994797& $-$& 0.077144\,i\\ [-8pt] \\
       147 & 0.400490& $-$& 0.286350\,i& 0.419276& $-$& 0.254325\,i& 
         1.001116& $-$& 0.075706\,i\\ [-8pt] \\
       148 & $-$0.091879& $+$& 0.669230\,i& $-$0.143660& $+$& 0.661309\,i& 
         0.995194& $-$& 0.077257\,i\\ [-8pt] \\
       149 & 0.946786& $-$& 0.263649\,i& 0.962119& $-$& 0.191329\,i& 
         0.999048& $-$& 0.075356\,i\\ [-8pt] \\
       150 & $-$0.024553& $-$& 0.058313\,i& $-$0.021860& $-$& 0.059113\,i
         & 1.002906& $-$& 0.044484\,i\\ [-8pt] \\
     \bottomrule \end{tabular}   
     \bigskip
    \caption{Exact and asymptotic values for $X=\Sigma (2,3,17)$}
  \end{table}

There is numerical evidence for Conjecture~\ref{thm:32} for certain Seifert
fibered manifolds\footnote{This has been superseded by the analytical work
in~\cite{Roz}, as explained in Remark~\ref{thm:26}.}~\cite{FG} and for some
hyperbolic examples as well~\cite{Ha},\cite{KSV}.  Table~1 compares
some explicit values of~$F_k(X)$ and the approximation given by the
asymptotic expansion in case $X$~is the Brieskorn homology sphere $X=\Sigma
(2,3,17)$.  The ratio of~$F_k(X)$ to the predicted asymptotic value is close
to one.  We include this chart partly as an antidote to the abstractions we
recounted from the formal properties of the path integral: the quantum
Chern-Simons invariants are calculable!  Indeed, these calculations were some
of the first concrete evidence in a mathematical context that path integrals
work as physicists claim.  Yet this most basic statement about the path
integral---the stationary phase asymptotic expansion---remains unproved in
this topological example.\footnote{See~\cite{CG} for some recent discussion.} 
 
There is much more to be said about the quantum Chern-Simons invariants in
mathematics.  For example, the deeper terms in the asymptotic expansion are
also invariants of 3-manifolds (of ``finite type''; see~\cite{Le} for a
survey) and links~\cite{B-N} .  The latter were invented by
Vassiliev~\cite{Va} in a different context.  There have been works too
numerous to mention which study the mathematical structure of these
invariants, including analytic properties of generating functions,
number-theoretical properties, asymptotics of knot invariants, relations to
hyperbolic geometry, etc.  Probably the most innovative proposal involving
Chern-Simons invariants is that they may be used in a model---both
theoretical and practical---for quantum computation~\cite{DFNSS}.  The
Chern-Simons construction appears often in theoretical physics.  One of the
simplest uses is as a mass term in 3-dimensional abelian gauge theory.  The
Chern-Simons invariant of the Levi-Civita connection also appears in quantum
theories of gravity~\cite{W2} in 3~dimensions.  Higher dimensional and higher
degree analogs of Chern-Simons appear in supergravity and string theory.  The
quantum Chern-Simons invariants are conjectured to have a relation to
Gromov-Witten invariants~\cite{GV}.  From the beginning they have been
closely connected with 2-dimensional conformal field theory.  And on and on.
 
Chern-Simons theory is one of many fronts over the last 25~years in a
vigorous interaction between geometry and theories of quantum fields and
quantum gravity.  It has been a very fruitful period.  The depth and variety
of mathematics which has been applied to problems in physics is astounding.
A new generation of theoretical physicists has a vast array mathematical
tools at its fingertips, and these are being used to investigate many
physical models.  While the mathematical community is justifiably excited
about new frontiers for applications of mathematics, such as biology, we
should also keep in mind the impressive successes that contemporary
applications of mathematics to fundamental physics enjoy.

What is particularly appealing for mathematicians in this interaction is the
impact in the other direction: theoretical physics has opened up many new
avenues of investigation in mathematics.  Quantum field theory and string
theory have made specific and concrete predictions which have motivated much
mathematical work over this period.  In part this is due to a temporary gulf
between theory and physical experiments: in place of the traditional
interplay between theorists and experimentalists has been an engagement with
mathematics. Mathematicians have tested predictions of quantum field theory
and string theory with independent mathematical techniques, in that way
serving as experimentalists.\footnote{As I write the LHC in Geneva is about
to turn on.  There is a possibility that new data from that large machine
will turn practitioners of modern mathematical physics back to a more
traditional connection with experiment.} But the ramifications for
mathematics go much deeper.  There are new connections between areas of
mathematics and new lines of research in existing areas.  It is a very
exciting time.
 
However, it is hardly the moment to declare victory.  What I find most
exciting are the possibilities for the future, as I believe we have only just
begun to absorb what this physics has to offer.  Let me be more specific.
One of the areas of success is \emph{topological} aspects of quantum field
theory, which includes the topological quantum field theories treated above.
The axiomatics touched upon in~\S\ref{sec:1} encode a part of quantum field
theory which applies beyond purely topological theories.  Indeed, it is an
old idea of 't Hooft\cite{'t} that much can be learned about physical
theories by exploring the topological information in the incarnations of a
theory at different scales, and this idea has been valuable in many contexts.
As we have seen, contemporary mathematical research in topological quantum
field theory is yielding new structure theorems and interesting applications.
Another exciting current area of math-physics interaction relates the
geometric Langlands program to four-dimensional gauge theories~\cite{KW}.
Lessons learned from 3-tier TQFTs, as represented in~\eqref{eq:10}, have had
a large impact here.\footnote{In this $n=4$~dimensional topological field
theory there is a category assigned to each Riemann surface.  That category
is argued to be one, or at least closely related to one, which enters the
geometric Langlands program.  There is a basic symmetry in 4-dimensional
gauge theory, S-duality, which relates two different theories, and on a
Riemann surface this is meant to induce the conjectured geometric Langlands
correspondence between two different categories attached to the surface.}
Another broad area of success over the past 25~years has been in
\emph{conformal} field theory, particularly in two dimensions.  There is a
similar geometric axiomatization~\cite{Se2} as well as more algebraic
approaches~\cite{FBZ}, \cite{BeDr} of which the three cited references are
only the tip of the iceberg.  Representation theory---most recently geometric
parts of representation theory---has been one part of mathematics intimately
connected with conformal field theory.  Also, mirror symmetry in
2-dimensional conformal field theory and string theory~\cite{HKKPTVVZ} has
had many consequences for algebraic geometry.
 
The common feature shared by topological and conformal theories is the
absence of \emph{scale}.  Conformal invariance is precisely the statement
that a theory looks the same at all scales.  And topology is usually studied
without introducing scale in the first place.  Yet scale is one of the most
basic concepts in all of physics, in particular in quantum field theory.
Even if one begins with a theory which is conformally invariant at the
classical level---and many important examples, such as Yang-Mills theory in
4-dimensions and $\sigma $-models in two dimensions, are---the process of
\emph{regularization} alluded to in~\S\ref{sec:5} introduces a scale into the
theory.  The dependence of the quantum theory on scale, encoded in the
\emph{renormalization group}~\cite{Wi}, \cite{Po} is a basic part of
quantum field theory which guides every practicing physicist in the field.
There is a body of mathematical work on quantum field theory (see~\cite{GJ}
for one account) which focuses on foundational analytic aspects which include
scale, but the newer investigations into topological and conformal features
have been largely disjoint.  The asymptotic expression~\eqref{eq:58} belongs
to the part of quantum field theory which depends on scale: it is
derived~\cite{W1} through a geometric version of regularization together with
standard manipulations of an integral.  Its unrelieved conjectural status is
only one of many signals that we need to turn more attention to the
scale-dependent aspects of quantum field theory.  I emphasize that this is
important even in theories which are scale-independent.  In particular, new
ideas are needed to integrate the scale-dependent and scale-independent
aspects.  
 
We illustrate with an episode from 4-dimensional topology and field theory.
In~1988, Witten~\cite{W3} introduced a topological twisting of a
supersymmetric gauge theory which encodes the Donaldson invariants~\cite{Do}
of 4-manifolds.  The underlying quantum field theory is not topological, but
certain correlation functions in a twisted version on oriented Riemannian
4-manifolds compute Donaldson invariants.  This new context for these
invariants of smooth manifolds bore fruit six years later when Seiberg and
Witten~\cite{SW} described the long-distance physics of the underlying
quantum field theory.  Their work is very much a part of scale-dependent
physics: the long-distance approximation is a field theory with different
fields than the fundamental short-distance theory.  Seiberg and Witten use
not only standard perturbation theory and renormalization group ideas, but
also the supersymmetry which severely constrains the form of the
long-distance physics.  As the Donaldson invariants are topological, so
independent of scale, they have an expression in both the long-distance and
short-distance theories.  In the short-distance theory one gets out
Donaldson's definition, which is Witten's original work, but the
long-distance approximation gives a new expression in terms of new
equations---what are now known as the Seiberg-Witten equations.  Geometers
immediately ran with these new equations, which are simpler than the
instanton equations used by Donaldson, and for many geometric purposes are
all one needs.  But it is important to realize that they are part of a
scale-dependent story in quantum field theory.  There has been considerable
effort to prove the conjectured equivalence between the two sets of
4-manifold invariants, and very recently it was finally proved by Feehan and
Leness~\cite{FL} for 4-manifolds of simple type.  (All known 4-manifolds are
of simple type.)  Here again we see that the \emph{discovery} of
scale-independent mathematics rests on a deep understanding of
scale-dependence.
 
The recent solution of the Poincar\'e conjecture by Perelman~\cite{P} offers
some insight.  The statement of the theorem---a simply connected 3-manifold
is homeomorphic to the 3-sphere---is very much in the realm of
scale-invariant mathematics: no scale enters at all.  Hamilton's Ricci
flow~\cite{H} introduces scale on a smooth 3-manifold in the form of a
Riemannian metric.  That metric may exhibit irregularities at arbitrary small
distances.  The Ricci flow smooths out the small-scale fluctuations in the
metric and, as time evolves, shifts focus to larger and larger scales.
Topological consequences, such as Poincar\'e, are deduced in the infinite
time limit when one approaches scale-independence.  Of course, this is a
poetic rendition of very intricate mathematics, especially as the flow may be
interrupted by singularities.  Still, it serves as a prototype for
scale-dependence in quantum field theories, and for good reason: Ricci flow
is the renormalization group flow for a particular quantum field theory, the
2-dimensional $\sigma $-model~\cite{Fr}.
 
The past 25~years have seen much interaction between physics and algebraic,
topological, and geometrical ideas.  The addition of analytic ideas relevant
to the study of scale will enrich this area in the years ahead.

 \appendix
  \section*{appendix: the chern-simons-weil theory of connections}\label{sec:4}

\setcounter{section}{1}
\setcounter{equation}{0}

Let $G$~be a Lie group with finitely many components.  Let $\mathfrak{g}$~be
its Lie algebra and $\pi \:P\to M$ a principal $G$-bundle.  Each
element~$g\in G$ acts on~$P$ by a diffeomorphism $R_g\:P\to P$.  Let $\theta
\in \Omega ^1_G(\mathfrak{g})$ be the left-invariant Maurer-Cartan form; it
transfers to any right $G$-torsor.\footnote{A \emph{torsor} for~$G$ is a
manifold on which $G$~acts simply transitively.  The fibers of a principal
bundle are right $G$-torsors.}  Denote by $i_m\:P_m \to P$ the inclusion of
the fiber at~$m\in M$.

        \begin{definition}[]\label{thm:11}
 A \emph{connection} on~$P$ is a 1-form $\Theta \in \Omega
^1_P(\mathfrak{g})$ which satisfies 
  \begin{equation}\label{eq:23}
     R_g^*\Theta =\Ad_{g\inv }\Theta ,\qquad i_m^*\Theta =\theta 
  \end{equation}
for all~$m\in M$. 
        \end{definition}

\noindent 
 The \emph{curvature} of~$\Theta $ is the $\mathfrak{g}$-valued 2-form 
  \begin{equation}\label{eq:24}
     \Omega =d\Theta + \frac 12[\Theta \wedge \Theta ]. 
  \end{equation}
It satisfies the linear equations 
  \begin{equation*}
     R_g^*\Omega =\Ad_{g\inv }\Omega,\qquad i_M^*\Omega =0,
  \end{equation*}
so lives in the linear space~$\Omega ^2_M(\mathfrak{g}\mstrut _{P})$ of
2-forms on~$M$ with values in the \emph{adjoint bundle} of Lie algebras
associated to~$P$.  Differentiating~\eqref{eq:24} we obtain the Bianchi
identity $d\Omega + [\Theta \wedge \Omega ]=0$.  To work out these equations
the reader will use the Jacobi identity $[[\Theta \wedge \Theta ]\wedge
\Theta ]=0$ and the structure equation $d\theta +\frac 12[\theta \wedge
\theta ]=0$.

The Chern-Weil construction is as follows.  Fix 
  \begin{equation}\label{eq:26}
     \langle \quad \rangle\:\mathfrak{g}^{\otimes p}\longrightarrow \RR 
  \end{equation}
which is $G$-invariant and symmetric. 

        \begin{proposition}[]\label{thm:12}
 $\omega (\Theta )=\langle \Omega \wedge \dots \wedge \Omega   \rangle$ is a
\emph{closed} $2p$-form on~$M$ 
        \end{proposition}

\noindent
 The proof is immediate: differentiate~$\omega (\Theta )$ and use Bianchi and
the $\mathfrak{g}$-invariance.    

        \begin{example}[]\label{thm:13}
 Let $\Sigma $~be an oriented Riemannian 2-manifold, $\pi \:P\to \Sigma $ its
$SO_2$-bundle of oriented orthonormal frames, and $\Theta $~the Levi-Civita
connection.  Then $\langle\Omega\rangle =K\,\vol\mstrut _\Sigma $, where
$\vol$~is the Riemannian volume 2-form.  In this case~$p=1$ and we define the
linear map~\eqref{eq:26} to identify the Lie algebra of~$SO_2$ with~$\RR$.
(We need not assume $M$~is oriented: the Chern-Weil construction for the
$O_2$-bundle of frames produces the density~$K\,d\mu \mstrut _\Sigma $ which
appears in~\eqref{eq:20}.)

        \end{example}

The de Rham theorem tells that $\omega (\Theta )$~determines a cohomology
class on~$M$.  It follows from~\eqref{eq:28} below that this class is
independent of the connection so is a topological invariant of the principal
bundle~$\pi \:P\to M$.  More important to us here is the geometry which
follows from the fact that the pullback~$\pi ^*\omega (\Theta )$ is exact.
The construction of Chern and Simons~\cite{CSi} which we now recount produces
a canonical antiderivative~\eqref{eq:29}, following their philosophy that
``the manner in which a closed form which is zero in cohomology actually
becomes exact contains geometric information'' (quoted
in~\cite{DGMS}).\footnote{This is analogous to the rise of category number in
geometry: the manner in which two sets are isomorphic contains geometric
information.}

As \eqref{eq:23}~are affine equations the space of solutions~$\AP\subset
\Omega ^1_P(\mathfrak{g})$ is also affine.  It is the space of all
connections on~$P$.  There is a universal connection~$\Tuniv$ on $\AP\times
P\to\AP\times M$ characterized by
  \begin{equation*}
     \Tuniv\res{\{\Theta \}\times P} = \Theta ,\qquad \Tuniv\res{\AP\times
     \{p\}}=0, \qquad \Theta \in \AP,\;p\in P. 
  \end{equation*}
Let $\Ouniv$~denote its curvature and $\omega (\Tuniv)$~its Chern-Weil form.
Now any two connections~$\Theta _0,\Theta _1\in \AP$ are the endpoints of a
line segment $\Delta ^1\to \AP$, an affine map with domain the standard
1-simplex.  Define the \emph{Chern-Simons form}
  \begin{equation}\label{eq:27}
     \alpha (\Theta _0,\Theta _1) = \int_{\Delta ^1}\omega (\Tuniv)\quad\in
     \Omega ^{2p-1}_M. 
  \end{equation}
Stokes' formula implies 
  \begin{equation}\label{eq:28}
     d\alpha (\Theta _0,\Theta _1) = \omega (\Theta _1) - \omega (\Theta
     _0). 
  \end{equation}
The affine structure of~$\AP$ allows us to continue to higher simplices:
three connections determine a 2-simplex and, by integration, a $(2p-2)$-form
whose differential relates the three Chern-Simons forms on the boundary,
etc.   
 
Our interest is to define a geometric invariant of a single connection.  For
this we work with the pullback bundle $\pi ^*P\to P$.  It has a canonical
trivialization $\Delta \:P\to \pi ^*P$: under the identification of $\pi ^*P$
with the fiber product~$P\times \mstrut _{M}P$ the section~$\Delta $ is the
diagonal map.  This section defines a trivial connection~$\TD\in \ApP$
characterized by~$\Delta ^*\TD=0$.  Pullback of connections defines an affine
embedding $\pi ^*\:\AP\hookrightarrow \ApP$.  For $\Theta \in \AP$ define the
\emph{Chern-Simons form}
  \begin{equation}\label{eq:29}
     \alpha (\Theta) = \alpha (\TD,\pi ^*\Theta )\quad\in \Omega ^{2p-1}_P. 
  \end{equation}
Then $\alpha (\Theta )$~is a differential form on~$P$.  From~\eqref{eq:28} we
deduce
  \begin{equation*}
     d\alpha (\Theta )= \pi ^*\omega (\Theta )\;\in \Omega ^{2p}_P. 
  \end{equation*}
The construction is functorial.  For example, applied to the unique
connection on the principal bundle $G\to\pt$ we obtain a closed form $\alpha
(\theta )\in \Omega ^{2p-1}_G$.  The functoriality, together
with~\eqref{eq:23}, implies the second of the equations
  \begin{equation*}
    R_g^*\alpha (\Theta )=\alpha (\Theta ),\qquad i_m^*\alpha (\Theta )=\alpha
    (\theta ).
  \end{equation*}
 
Both~\eqref{eq:27} and~\eqref{eq:29} are termed `Chern-Simons forms'; the
former lives on the base~$M$ and depends on two connections, the latter lives
on the total space~$P$ and depends on a single connection.

        \begin{example}[]\label{thm:14}
 In the situation of Example~\ref{thm:13} let $C\subset \Sigma $ be a smooth
oriented curve.  There is a canonical framing of~$\Sigma $ along~$C$---the
oriented tangent vector is the first basis vector---and so a canonical lift
of~$C$ to the frame bundle~$P$.  The Chern-Simons 1-form~$\alpha (\Theta )$
then pulls down to a 1-form on~$C$ which is $\kappa \mstrut _C\vol_C$;
compare~\eqref{eq:21}.
        \end{example}

The \emph{Chern-Simons invariant} generalizes the total geodesic
curvature~\eqref{eq:21}.  Suppose $\Theta $~is a connection on a principal
$G$-bundle $\pi \:P\to M$ with $M$~a closed oriented $(2p-1)$-manifold.  Let
$\langle \quad \rangle$~be an invariant polynomial as in~\eqref{eq:26} and
assume it is normalized so that the form $\alpha (\theta )\in \Omega
^{2p-1}_G$ has integral periods.  Assume first that $\pi \:P\to M$ is
trivializable.  Then if $s\:M\to P$ is any section of~$\pi $, the
Chern-Simons invariant
  \begin{equation*}
     S(\Theta )= \int_{M}s^*\alpha (\Theta )\pmod1 
  \end{equation*}
is independent of~$s$.  Note that the invariant~$S(\Theta )$ lives
in~$\RR/\ZZ$. 
 
For a general (nontrivializable) bundle $\pi \:P\to M$ we need a more
intricate construction~\cite{F3}, \cite{DW}.  Let $EG\to BG$ be a universal
bundle built from smooth (Hilbert) manifolds, $\Tun$~a connection on~$EG$,
and $\gamma \:P\to EG$ a $G$-equivariant (classifying) map with quotient
$\bg\:M\to BG$.  Since the odd homology of~$BG$ is torsion, for some positive
integer~$N$ there is a smooth $2p$-chain~$W$ in~$BG$ with $\partial W=N\cdot
\bg([M])$.  Then $[\frac 1N W]\in H_{2p}(BG;\RR/\ZZ)$.  The Chern-Simons
invariant depends on a choice of cohomology class $\lambda \in
H^{2p}(BG;\ZZ)$ whose image in~$H^{2p}(BG;\RR)$ is represented by the
universal Chern-Weil form $\omega (\Tun)$.  The integrality assumption
on~$\langle \quad   \rangle$ guarantees the existence of~$\lambda $.
For a connection~$\Theta $ on~$P$ set 
  \begin{equation*}
     S(\Theta ) = \frac 1N\int_{W}\omega (\Tun) + \int_{M}\alpha (\gamma
     ^*\Tun,\Theta ) + \lambda \bigl([\frac 1N W] \bigr)\quad \in \RR/\ZZ. 
  \end{equation*}
One can check that the right hand side is independent of the choices
of~$\gamma ,W$ and defines a smooth function $S\:\AP\to\RR/\ZZ $.

\bigskip\bigskip

\newcommand{\etalchar}[1]{$^{#1}$}
\providecommand{\bysame}{\leavevmode\hbox to3em{\hrulefill}\thinspace}
\providecommand{\MR}{\relax\ifhmode\unskip\space\fi MR }
\providecommand{\MRhref}[2]{%
  \href{http://www.ams.org/mathscinet-getitem?mr=#1}{#2}
}
\providecommand{\href}[2]{#2}

\end{document}